\documentclass[11pt]{article}

\usepackage{geometry}
\usepackage{latexsym}
\usepackage{mathrsfs}
\usepackage{amsfonts}
\usepackage{amssymb}
\usepackage{subfigure}    
\usepackage{graphicx}
\usepackage{epstopdf}
\usepackage{float}
\usepackage{multirow}
\usepackage{bm}
\usepackage{amsmath}
\usepackage{amsthm}
\usepackage{color}
\usepackage[subnum]{cases}
\usepackage{rotating}
\usepackage{enumitem}
\usepackage{accents}
\usepackage{algorithm}
\usepackage{algorithmic}
\usepackage[title]{appendix}
\usepackage{mathtools}
\usepackage{caption}
\usepackage{tcolorbox}

\usepackage{url}
\usepackage[pdftex,
plainpages=false,
bookmarks,
bookmarksnumbered,
colorlinks=true,
linkcolor=blue,
citecolor=blue,
urlcolor=blue,
filecolor=black
]
{hyperref}

\makeatletter
\@addtoreset{equation}{section}
\makeatother

\def \beginproof{\par\noindent {\bf Proof.}\ }
\def \endproof{\hskip .5cm $\Box$ \vskip .5cm}

\makeatletter
\def \wideubar{\underaccent{{\cc@style\underline{\mskip15mu}}}}
\def \widebar{\accentset{{\cc@style\underline{\mskip11mu}}}}
\makeatother

\graphicspath{{./images/log/}{./images/l1l2/}}

\begin{document}

\newtheorem{property}{Property}[section]
\newtheorem{proposition}{Proposition}[section]
\newtheorem{append}{Appendix}[section]
\newtheorem{definition}{Definition}[section]
\newtheorem{lemma}{Lemma}[section]
\newtheorem{corollary}{Corollary}[section]
\newtheorem{theorem}{Theorem}[section]
\newtheorem{remark}{Remark}[section]
\newtheorem{problem}{Problem}[section]
\newtheorem{example}{Example}[section]
\newtheorem{assumption}{Assumption}[section]
\renewcommand*{\theassumption}{\Alph{assumption}}

\title{Proximal Gradient Method with Extrapolation and Line Search for a Class of Nonconvex and Nonsmooth Problems}
\author{Lei Yang\footnotemark[1]}
\date{}
\maketitle

\renewcommand{\thefootnote}{\fnsymbol{footnote}}
\footnotetext[1]{Department of Applied Mathematics, The Hong Kong Polytechnic University, Hung Hom, Kowloon, Hong Kong, China ({\tt yanglei.math@gmail.com}).}
\renewcommand{\thefootnote}{\arabic{footnote}}

\begin{abstract}
In this paper, we consider a class of possibly nonconvex, nonsmooth, and non-Lipschitz optimization problems arising in many contemporary applications such as machine learning, variable selection, and image processing. To solve this class of problems, we propose a proximal gradient method with extrapolation and line search (PGels). This method is developed based on a special potential function and successfully incorporates both \textit{extrapolation} and \textit{non-monotone line search}, which are two simple and efficient acceleration techniques for the proximal gradient method. Thanks to the line search, this method allows more flexibilities in choosing the extrapolation parameters and updates them adaptively at each iteration if a certain line search criterion is not satisfied. Moreover, with proper choices of parameters, our PGels reduces to many existing algorithms. We also show that, under some mild conditions, our line search criterion is well-defined and any cluster point of the sequence generated by PGels is a stationary point of our problem. In addition, by making assumptions on the Kurdyka-{\L}ojasiewicz exponent of the objective in our problem, we further analyze the local convergence rate of two special cases of the PGels, including the widely-used non-monotone proximal gradient method as one case. Finally, we conduct some preliminary numerical experiments for solving the $\ell_1$ regularized logistic regression problem and the $\ell_{1\text{-}2}$ regularized least squares problem. The numerical results show the promising performance of the PGels and illustrate the potential advantage of combining two acceleration techniques.

\vspace{5mm}

\end{abstract}


\section{Introduction}\label{secintro}

In this paper, we consider the following composite optimization problem:
\begin{eqnarray}\label{model}
\min\limits_{\bm{x}\in\mathbb{R}^n}~F(\bm{x}) := f(\bm{x}) + P(\bm{x}),
\end{eqnarray}
where $f: \mathbb{R}^n \rightarrow \mathbb{R}$ and $P: \mathbb{R}^n \rightarrow \mathbb{R}\cup\{+\infty\}$ satisfy Assumption \ref{assumfun}.

\begin{assumption}\label{assumfun}
\indent
\begin{itemize}
\item[{\rm (i)}] $f: \mathbb{R}^n \rightarrow \mathbb{R}$ is a continuously differentiable (possibly nonconvex) function with Lipschitz continuous gradient, i.e., there exists a Lipschitz constant $L_f > 0$ such that
                \begin{eqnarray*}
                \|\nabla f(\bm{x}) - \nabla f(\bm{y})\| \leq L_f \|\bm{x} - \bm{y}\|, \qquad \forall \bm{x}, \,\bm{y} \in \mathbb{R}^n.
                \end{eqnarray*}

\item[{\rm (ii)}] $P: \mathbb{R}^n \rightarrow \mathbb{R}\cup\{+\infty\}$ is a proper closed (possibly nonconvex, nonsmooth, and non-Lipschitz) function; it is bounded below and continuous on its domain. Moreover, the proximal mapping of $\nu P$ is easy to compute for all $\nu>0$ (see next section for notation and definitions).

\item[{\rm (iii)}] $F$ is level-bounded.
\end{itemize}
\end{assumption}

Problem \eqref{model} arises in many contemporary applications such as machine learning \cite{cs2017optimization,snw2012optimization}, variable selection \cite{fl2001variable,hhjm2008asymptotic,kf2000asymptotics,t1996regression,z2010nearly} and image processing \cite{cp2016introduction,nnzc2008efficient}. In general, $f$ is a loss or fitting function used for measuring the deviation of a solution from the observations. Two commonly used loss functions are the least squares loss function  $f(\bm{x})=\frac{1}{2}\|A\bm{x}-\bm{b}\|^2$ and the logistic loss function $f(\bm{x})=\sum^{m}_{i=1} \log(1+\exp(-b_i(\bm{a}_i^{\top}\bm{x})))$,
where $A = [\bm{a}_1,\cdots,\bm{a}_m]^{\top} \in \mathbb{R}^{m \times n}$ is a data matrix and $\bm{b}\in\mathbb{R}^m$ is an observed vector. One can easily verify that these two loss functions satisfy Assumption \ref{assumfun}(i). On the other hand, $P$ is a regularizer used for inducing certain structure in the solution. For example, $P$ can be the indicator function for a certain set such as $\mathcal{X}=\{\bm{x}\in\mathbb{R}^n : \bm{x} \geq 0\}$ and $\mathcal{X}=\{\bm{x}\in\mathbb{R}^n : \sum^n_{i=1} x_i = 1, \,\bm{x}\geq0\}$; the former choice restricts the elements of the solution to be nonnegative and the latter choice restricts the solution in a simplex. We can also choose $P$ to be a certain sparsity-inducing regularizer such as $\lambda\|\bm{x}\|_p^p$ for $0<p\leq1$ \cite{hhjm2008asymptotic,kf2000asymptotics,t1996regression}, $\lambda\sum^n_{i=1}\log(1+\alpha|x_i|)$ for $\alpha>0$ \cite{nnzc2008efficient} and $\lambda(\|\bm{x}\|_1-\|\bm{x}\|)$ \cite{ylhx2015Minimization}, where $\lambda>0$ is a regularization parameter. Note that all the aforementioned examples of $P$ as well as many other widely-used regularizers (see \cite{apx2017difference,c2012smoothing} and references therein for more regularizers) satisfy Assumption \ref{assumfun}(ii). Finally, we would like to point out that Assumption \ref{assumfun}(iii) is also satisfied by many choices of $f$ and $P$ in practice; see, for example, problems \eqref{l1logmodel} and \eqref{modell12} in our numerical part. More examples of problem \eqref{model} can be found in \cite{cp2016introduction,cs2017optimization,snw2012optimization} and references therein.

Due to the importance and the popularity of problem \eqref{model}, tremendous attempts have been made to solve it efficiently, especially when the problem involves a large number of variables. One popular class of methods for solving problem \eqref{model} is the first-order method due to the computational simplicity and good convergence properties. Among various first-order methods, the proximal gradient (PG) method \cite{fm1981a,lm1979splitting} (also known as the forward-backward splitting algorithm \cite{cp2011proximal} \textit{or} the iterative shrinkage-thresholding algorithm \cite{bt2009a}) is arguably the most fundamental one, whose basic iterative step reads as follows:
\begin{eqnarray}\label{pgscheme}
\bm{x}^{k+1} \in \mathop{\mathrm{Argmin}}\limits_{\bm{x}} \left\{\langle \nabla f(\bm{x}^k), \,\bm{x}\rangle + \frac{\mu}{2}\|\bm{x} - \bm{x}^k\|^2 + P(\bm{x})\right\},
\end{eqnarray}
where $\mu > 0$ is a constant depending on the Lipschitz constant $L_f$ of $\nabla f$. However, the PG method can be slow in practice; see, for example, \cite{bt2009a,oc2015adaptive,wcp2017linear}. Therefore, a large amount of research has been conducted to accelerate PG for solving problem \eqref{model}.

One simple and widely studied strategy is to perform \textit{extrapolation} in the spirit of Nesterov's acceleration techniques \cite{n1983a,n2013introductory}, whose key idea is to make use of historical information at each iteration. A typical scheme of the proximal gradient method with extrapolation (PGe) for solving problem \eqref{model} is
\begin{eqnarray}\label{pgescheme}
\left\{\begin{aligned}
&\bm{y}^k = \bm{x}^k + \beta_k (\bm{x}^k - \bm{x}^{k-1}),  \\
&\bm{x}^{k+1} \in \mathop{\mathrm{Argmin}}\limits_{\bm{x}} \left\{\langle \nabla f(\bm{y}^k), \,\bm{x}\rangle + \frac{\mu}{2}\|\bm{x} - \bm{y}^k\|^2 + P(\bm{x})\right\},
\end{aligned}\right.
\end{eqnarray}
where $\beta_k$ is the extrapolation parameter satisfying certain conditions and $\mu > 0$ is a constant depending on $L_f$. One representative algorithm that takes the form of \eqref{pgescheme} with proper choices of $\{\beta_k\}$ is the fast iterative shrinkage-thresholding algorithm (FISTA) \cite{bt2009a}, which extends Nesterov's classical accelerated gradient method \cite{n1983a} to a general composite setting
and exhibits a faster convergence rate of $O(1/k^2)$ in terms of objective values (see \cite{bt2009a,n1983a} for more details). This motivates the study of the PGe and its variants for solving problem \eqref{model} under different scenarios; see, for example, \cite{bt2009fast,bcg2011templates,gl2016accelerated,ll2015accelerated,ocbp2014ipiano,oc2015adaptive,t2008on,t2010approximation,wcp2017linear,wcp2017a,xy2013a,xy2015block,xy2017globally}. It is worth noting that, when $f$ and $P$ are convex, many existing PGe and its variants \cite{bt2009fast,bt2009a,bcg2011templates,n2013gradient,oc2015adaptive,t2008on,t2010approximation} choose the extrapolation parameters $\{\beta_k\}$ (explicitly or implicitly) based on the following updating scheme
\begin{eqnarray}\label{exparupdate}
\left\{\begin{aligned}
&\beta_k = (t_{k-1} - 1)/t_k,  \\
&t_{k+1} = \frac{1+\sqrt{1+4t_k^2}}{2},
\end{aligned}\right.
\quad \mathrm{with} \quad t_{-1}=t_0=1,
\end{eqnarray}
which originated from Nesterov's seminal works \cite{n1983a,n2013introductory} and was shown to be ``optimal" \cite{n2013introductory}. However, for the nonconvex case, the ``optimal" choices of $\{\beta_k\}$ are still not clear. Although the convergence of the PGe and its variants can be guaranteed in theory for some classes of nonconvex problems under certain conditions on $\{\beta_k\}$ (see, for example, \cite{gl2016accelerated,ll2015accelerated,ocbp2014ipiano,wcp2017linear,wcp2017a,xy2013a,xy2015block,xy2017globally}), the choices of $\{\beta_k\}$ are relatively restrictive and may not work well for achieving acceleration.

Another efficient strategy for accelerating the PG method is to employ a \textit{non-monotone line search technique} to adaptively find a proper $\mu$ in the scheme \eqref{pgscheme} at each iteration. The non-monotone line search technique dates back to the non-monotone Newton's method proposed by Grippo et al. \cite{gll1986nonmonotone} and has been applied to many algorithms with good empirical performances; see, for example, \cite{bmr2000nonmonotone,d2002nonmonotone,ypc2017a,zw2004nonmonotone}. Based on this technique, Wright et al. \cite{wnf2009sparse} recently proposed an efficient method (called SpaRSA) to solve problem \eqref{model}, whose iteration is roughly given as follows: Choose $\tau>1$, $c>0$ and an integer $N\geq0$. Then, at the $k$-th iteration, choose $\mu_k^0>0$ and find the smallest nonnegative integer $j_k$ such that
\begin{eqnarray*}
\left\{\begin{aligned}
&\bm{u} \in \mathop{\mathrm{Argmin}}\limits_{\bm{x}} \left\{\langle \nabla f(\bm{x}^k), \,\bm{x}\rangle + \frac{\tau^{j_k}\mu_k^0}{2}\|\bm{x} - \bm{x}^k\|^2 + P(\bm{x})\right\}, \\
&F(\bm{u}) - \max\limits_{[k-N]_{+}\leq i\leq k}F(\bm{x}^{i}) \leq -\frac{c}{2}\|\bm{u}-\bm{x}^{k}\|^2.
\end{aligned}\right.
\end{eqnarray*}
This method is essentially the non-monotone proximal gradient (NPG) method, namely, the proximal gradient method with a non-monotone line search. Later on, the NPG method was extended for solving problem \eqref{model} under more general conditions and has been shown to exhibit promising numerical performances in many application problems (see, for example, \cite{clp2016penalty,gzlhy2013a,lpt2017a}). In view of the above, it is natural to raise a question:


\begin{tcolorbox}[colback=black!0!white]
\textit{Can we derive an efficient method for solving problem \eqref{model}, which takes advantage of both extrapolation and non-monotone line search?}
\end{tcolorbox}

In this paper, we initiate our study to address this question. We propose a method for solving problem \eqref{model} that successfully incorporates both extrapolation and non-monotone line search and allows more flexibilities in the choices of the extrapolation parameters $\{\beta_k\}$. We shall call our method the proximal gradient method with extrapolation and line search, denoted by PGels for short. This method is developed based on the following potential function (specially constructed for the objective $F$ in problem \eqref{model}):
\begin{eqnarray}\label{defpofun}
H_{\delta}(\bm{u}, \bm{v}, \mu) := F(\bm{u}) + \frac{\delta\mu}{4}\|\bm{u} - \bm{v}\|^2, \quad \forall \,\bm{u},\,\bm{v} \in \mathbb{R}^n, ~\mu > 0,
\end{eqnarray}
where $\delta \in [0, \,1)$ is a given nonnegative constant. Clearly, $H_{\delta}(\bm{u}, \bm{v}, \mu) \equiv F(\bm{u})$ if $\delta=0$. We will see in Section \ref{secalg} that this potential function is used to establish a new non-monotone line search criterion \eqref{lscond} when the extrapolation technique is applied. This allows more choices of $\beta_k$ at each iteration, and will adaptively update $\mu_k$ and $\beta_k$ at the same time if the line search criterion is not satisfied (see Algorithm \ref{alg_PGels} for more details). The convergence analysis of the PGels is also presented in Section \ref{secalg}. Specifically, under Assumption \ref{assumfun}, we show that our line search criterion \eqref{lscond} is well-defined and any cluster point of the sequence generated by the PGels is a stationary point of problem \eqref{model}. Moreover, since our PGels readily reduces to PG, PGe or NPG with proper choices of parameters (see Remark \ref{eqalgo}), then we actually obtain a unified convergence analysis for PG, PGe and NPG as a byproduct. In addition, in Section \ref{secconrate}, we further study the local convergence rate in terms of objective values for two special cases of the PGels (including the NPG as one case) under an additional assumption on the Kurdyka-{\L}ojasiewicz exponent of the objective $F$ in problem \eqref{model}. To the best of our knowledge, this is the first local convergence rate analysis of the NPG for solving problem \eqref{model}. Finally, we conduct some numerical experiments in Section \ref{secnum} to evaluate the performance of our method for solving the $\ell_1$ regularized logistic regression problem and the $\ell_{1\text{-}2}$ regularized least squares problem. Our computational results illustrate the efficiency of our method and show the potential advantage of combining extrapolation and non-monotone line search.

The rest of this paper is organized as follows. In Section \ref{secnot}, we present notation and preliminaries used in this paper. In Section \ref{secalg}, we describe the PGels for solving problem \eqref{model} and study its global subsequential convergence. The local convergence rate of two special cases of the PGels is analyzed in Section \ref{secconrate} and some numerical results are reported in Section \ref{secnum}. Finally, some concluding remarks are given in Section \ref{secconc}.

\section{Notation and preliminaries}\label{secnot}

In this paper, we present scalars, vectors and matrices in lower case letters, bold lower case letters and upper case letters, respectively. We also use $\mathbb{R}$, $\mathbb{R}^n$, $\mathbb{R}^n_+$ and $\mathbb{R}^{m\times n}$ to denote the set of real numbers, $n$-dimensional real vectors, $n$-dimensional real vectors with nonnegative entries and $m\times n$ real matrices, respectively. For a vector $\bm{x}\in\mathbb{R}^n$, $x_i$ denotes its $i$-th entry, $\|\bm{x}\|$ denotes its Euclidean norm, $\|\bm{x}\|_1$ denotes its $\ell_1$ norm defined by $\|\bm{x}\|_1:=\sum^n_{i=1}|x_i|$, $\|\bm{x}\|_p$ denotes its $\ell_p$-quasi-norm ($0 < p < 1$) defined by $\|\bm{x}\|_p:=\left( \sum_{i=1}^{n} |x_{i}|^p \right)^{\frac{1}{p}}$ and $\|\bm{x}\|_{\infty}$ denotes its $\ell_{\infty}$ norm given by the largest entry in magnitude. For a matrix $A\in\mathbb{R}^{m\times n}$, its spectral norm is denoted by $\|A\|$, which is the largest singular value of $A$.

For an extended-real-valued function $h: \mathbb{R}^{n} \rightarrow [-\infty,\infty]$, we say that it is \textit{proper} if $h(\bm{x}) > -\infty$ for all $\bm{x} \in \mathbb{R}^{n}$ and its domain ${\rm dom}\,h:=\{\bm{x} \in \mathbb{R}^{n} : h(\bm{x}) < \infty\}$ is nonempty. A proper function $h$ is said to be closed if it is lower semicontinuous. We also use the notation $\bm{y} \xrightarrow{h} \bm{x}$ to denote $\bm{y} \rightarrow \bm{x}$ and $h(\bm{y}) \rightarrow h(\bm{x})$. The basic \textit{subdifferential} (see \cite[Definition 8.3]{rw1998variational}) of $h$ at $\bm{x} \in \mathrm{dom}\,h$ used in this paper is
\begin{eqnarray*}
\partial h(\bm{x}) := \left\{ \bm{d} \in \mathbb{R}^{n}: \exists \,\bm{x}^k \xrightarrow{h} \bm{x}, ~\bm{d}^k \rightarrow \bm{d} ~~\mathrm{with}~\liminf\limits_{\bm{y} \rightarrow \bm{x}^k, \,\bm{y} \neq \bm{x}^k}\, \frac{h(\bm{y})-h(\bm{x}^k)-\langle \bm{d}^k, \bm{y}-\bm{x}^k\rangle}{\|\bm{y}-\bm{x}^k\|} \geq 0\ ~\forall k\right\}.
\end{eqnarray*}
It can be observed from the above definition that
\begin{eqnarray}\label{robust}
\left\{ \bm{d}\in\mathbb{R}^{n}: \exists \,\bm{x}^k \xrightarrow{h} \bm{x}, ~\bm{d}^k \rightarrow \bm{d} ~\mathrm{with}~\bm{d}^k \in \partial h(\bm{x}^k)~\mathrm{for}~\mathrm{each}~k \right\} \subseteq \partial h(\bm{x}).
\end{eqnarray}
When $h$ is continuously differentiable or convex, the above subdifferential coincides with the classical concept of derivative or convex subdifferential of $h$; see, for example, \cite[Exercise~8.8]{rw1998variational} and \cite[Proposition~8.12]{rw1998variational}. In addition, if $h$ has several groups of variables, we use $\partial_{\bm{x}_i} h$ (resp., $\nabla_{\bm{x}_i} h$) to denote the partial subdifferential (resp., gradient) of $h$ with respect to the group of variables $\bm{x}_i$.

For a proper closed function $h: \mathbb{R}^{n} \rightarrow (-\infty,\infty]$ and $\nu>0$, the proximal mapping of $\nu h$ at $\bm{y}\in\mathbb{R}^n$ is defined by
\begin{eqnarray*}
\mathrm{Prox}_{\nu h}(\bm{y}) := \mathop{\mathrm{Argmin}}\limits_{\bm{x}\in\mathbb{R}^n} \left\{h(\bm{x}) + \frac{1}{2\nu}\|\bm{x} - \bm{y}\|^2\right\}.
\end{eqnarray*}
Note that this operator is well defined for any $\bm{y}\in\mathbb{R}^n$ if $h$ is bounded below in $\mathbb{R}^n$. For a closed set $\mathcal{X}\subseteq\mathbb{R}^{n}$, its indicator function $\delta_{\mathcal{X}}$ is defined by
\begin{eqnarray*}
\delta_{\mathcal{X}}(\bm{x}) = \left\{
\begin{array}{ll}
0, &\quad\mathrm{if}~\bm{x}\in\mathcal{X}, \\
+\infty, &\quad\mathrm{otherwise}.
\end{array}\right.
\end{eqnarray*}
We also use $\mathrm{dist}(\bm{x}, \mathcal{X})$ to denote the distance from $\bm{x}$ to $\mathcal{X}$, i.e., $\mathrm{dist}(\bm{x}, \mathcal{X}) := \inf_{\bm{y}\in\mathcal{X}}\|\bm{x}-\bm{y}\|$.

For any local minimizer $\bar{\bm{x}}$ of \eqref{model}, it is known from \cite[Theorem~10.1]{rw1998variational} and \cite[Exercise~8.8(c)]{rw1998variational} that the following first-order necessary condition holds:
\begin{eqnarray}\label{optcond}
0 \in \nabla f(\bar{\bm{x}}) + \partial P(\bar{\bm{x}}),
\end{eqnarray}
where $\nabla f$ denotes the gradient of $f$. In this paper, we say that $\bm{x}^*$ is a \textit{stationary point} of \eqref{model} if $\bm{x}^*$ satisfies \eqref{optcond} in place of $\bar{\bm{x}}$.

We next recall the Kurdyka-{\L}ojasiewicz (KL) property (see \cite{ab2009on,abrs2010proximal,abs2013convergence,bdl2007the,bst2014proximal} for more details), which plays an important role in our analysis for the local convergence rate in Section \ref{secconrate}. For notational simplicity, let $\Xi_{\nu}$ ($\nu>0$) denote a class of concave functions $\varphi:[0,\nu) \rightarrow \mathbb{R}_{+}$ satisfying: (i) $\varphi(0)=0$; (ii) $\varphi$ is continuously differentiable on $(0,\nu)$ and continuous at $0$; (iii) $\varphi'(t)>0$ for all $t\in(0,\nu)$. Then, the KL property can be described as follows.

\begin{definition}[\textbf{KL property and KL function}]
Let $h: \mathbb{R}^n \rightarrow \mathbb{R} \cup \{+\infty\}$ be a proper closed function.
\begin{itemize}
\item[(i)] For $\tilde{\bm{x}}\in{\rm dom}\,\partial h:=\{\bm{x} \in \mathbb{R}^{n}: \partial h(\bm{x}) \neq \emptyset\}$, if there exist a $\nu\in(0, +\infty]$, a neighborhood $\mathcal{V}$ of $\tilde{\bm{x}}$ and a function $\varphi \in \Xi_{\nu}$ such that for all $\bm{x} \in \mathcal{V} \cap \{\bm{x}\in \mathbb{R}^{n} : h(\tilde{\bm{x}})<h(\bm{x})<h(\tilde{\bm{x}})+\nu\}$, it holds that
    \begin{eqnarray*}
    \varphi'(h(\bm{x})-h(\tilde{\bm{x}}))\,\mathrm{dist}(0, \,\partial h(\bm{x})) \geq 1,
    \end{eqnarray*}
    then $h$ is said to have the \textbf{Kurdyka-{\L}ojasiewicz (KL)} property at $\tilde{\bm{x}}$.

\item[(ii)] If $h$ satisfies the KL property at each point of ${\rm dom}\,\partial h$, then $h$ is called a KL function.
\end{itemize}
\end{definition}

A large number of functions such as proper closed semialgebraic functions satisfy the KL property \cite{abrs2010proximal,abs2013convergence}. Based on the above definition, we then introduce the KL exponent \cite{abrs2010proximal,lp2017calculus}.

\begin{definition}[\textbf{KL exponent}]
Suppose that $h: \mathbb{R}^n \rightarrow \mathbb{R} \cup \{+\infty\}$ is a proper closed function satisfying the KL property at $\tilde{\bm{x}}\in{\rm dom}\,\partial h$ with $\varphi(t)=a' t^{1-\theta}$ for some $a' > 0$ and $\theta\in[0, 1)$, i.e., there exist $a, \varepsilon, \nu > 0$ such that
\begin{eqnarray*}
\mathrm{dist}(0, \,\partial h(\bm{x})) \geq a \left(h(\bm{x}) - h(\tilde{\bm{x}})\right)^{\theta}
\end{eqnarray*}
whenever $\bm{x} \in {\rm dom}\,\partial h$, $\|\bm{x}-\tilde{\bm{x}}\| \leq \varepsilon$ and $h(\tilde{\bm{x}})<h(\bm{x})<h(\tilde{\bm{x}})+\nu$. Then, $h$ is said to have the KL property at $\tilde{\bm{x}}$ with an exponent $\theta$. If $h$ is a KL function and has the same exponent $\theta$ at any $\tilde{\bm{x}}\in{\rm dom}\,\partial h$, then $h$ is said to be a KL function with an exponent $\theta$.
\end{definition}

We also recall the uniformized KL property, which was established in \cite[Lemma 6]{bst2014proximal}.

\begin{proposition}[\textbf{Uniformized KL property}]\label{uniKL}
Suppose that $h: \mathbb{R}^n \rightarrow \mathbb{R} \cup \{+\infty\}$ is a proper closed function and $\Gamma$ is a compact set. If $h \equiv \zeta$ on $\Gamma$ for some constant $\zeta$ and satisfies the KL property at each point of $\Gamma$, then there exist $\varepsilon>0$, $\nu>0$ and $\varphi \in \Xi_{\nu}$ such that
\begin{eqnarray*}
\varphi'(h(\bm{x}) - \zeta)\,\mathrm{dist}(0, \,\partial h(\bm{x})) \geq 1
\end{eqnarray*}
for all $\bm{x} \in \{\bm{x}\in\mathbb{R}^{n}: \mathrm{dist}(\bm{x},\,\Gamma)<\varepsilon\} \cap \{\bm{x}\in \mathbb{R}^{n} : \zeta < h(\bm{x}) < \zeta + \nu\}$.
\end{proposition}

Finally, we give two supporting lemmas that will be used in our analysis.

\begin{lemma}[{\cite[Lemma 2.2]{lp2017calculus}}]\label{normineq1}
Let $\alpha > 0$. Then, for any $\bm{w} = (w_1, \cdots, w_n)^{\top}\in\mathbb{R}^n_{+}$, there exist $0 < c_1 \leq c_2$ such that $c_1\|\bm{w}\| \leq \left(w_1^{\alpha}+\cdots+w_n^{\alpha}\right)^{\frac{1}{\alpha}} \leq c_2 \|\bm{w}\|$.
\end{lemma}

\begin{lemma}[{\cite[Lemma 3.1]{lp2017calculus}}]\label{normineq2}
Let $1 < \alpha \leq 2$. Then, for any $\bm{u}, \bm{v}\in\mathbb{R}^n$, there exist $b_1 > 0$ and $0< b_2 < 1$ such that $\|\bm{u}+\bm{v}\|^{\alpha} \geq b_1 \|\bm{u}\|^{\alpha}-b_2\|\bm{v}\|^{\alpha}$.
\end{lemma}

\section{Proximal gradient method with extrapolation and line search}\label{secalg}

In this section, we present a proximal gradient method with extrapolation and line search (PGels) for solving problem \eqref{model}. This method is developed based on a specially constructed potential function defined in \eqref{defpofun}. The complete framework is presented in Algorithm \ref{alg_PGels}.

\begin{algorithm}[ht]
\caption{PGels for solving problem \eqref{model}}\label{alg_PGels}
\textbf{Input:} $\bm{x}^0\in\mathrm{dom}\,F$, $\tau>1$, $0\leq\delta<1$, $0<\eta<1$, $c>0$, $\mu_{\max}\geq\frac{L_f+2c}{1-\delta}\geq\mu_{\min}>0$, $\beta_{\max}\geq0$ and an integer $N\geq0$. Set $\bm{x}^{-1}=\bm{x}^0$, $\bar{\mu}_{-1}=1$ and $k=0$. \vspace{1mm}\\
\textbf{while} a termination criterion is not met, \textbf{do}
\begin{itemize}[leftmargin=2cm]
\item[\textbf{Step 1}.] Choose $\mu_k^0\in[\mu_{\min}, \,\mu_{\max}]$ and $\beta_k^0\in[0,\,\delta\beta_{\max}]$ arbitrarily. Set $\mu_k = \mu_k^0$ and $\beta_k = \beta_k^0$.
                        \begin{itemize}
                        \item[\textbf{(1a)}] Compute
                                             \begin{eqnarray}\label{extrastep}
                                             \bm{y}^k = \bm{x}^k + \beta_k (\bm{x}^k - \bm{x}^{k-1}).
                                             \end{eqnarray}

                        \item[\textbf{(1b)}] Solve the subproblem
                                             \begin{eqnarray}\label{subpro}
                                             \bm{u} \in \mathop{\mathrm{Argmin}}\limits_{\bm{x}} \left\{\langle \nabla f(\bm{y}^k), \,\bm{x} - \bm{y}^k\rangle + \frac{\mu_k}{2}\|\bm{x} - \bm{y}^k\|^2 + P(\bm{x})\right\}.
                                             \end{eqnarray}

                        \item[\textbf{(1c)}] If
                                             \begin{eqnarray}\label{lscond}
                                             \begin{aligned}
                                             H_{\delta}(\bm{u}, \bm{x}^k, \mu_k) - \max\limits_{[k-N]_{+}\leq i\leq k}H_{\delta}(\bm{x}^{i}, \bm{x}^{i-1}, \bar{\mu}_{i-1})\leq-\frac{c}{2}\|\bm{u}-\bm{x}^{k}\|^2
                                             \end{aligned}
                                             \end{eqnarray}
                                             is satisfied, then go to \textbf{Step 2}.

                        \item[\textbf{(1d)}] Set $\mu_k\leftarrow\min\{\tau\mu_k, \,\mu_{\max}\}$, $\beta_k\leftarrow\eta\beta_k$ and go to step \textbf{(1a)}.
                        \end{itemize}

\item [\textbf{Step 2}.] Set $\bm{x}^{k+1} \leftarrow \bm{u}$, $\bar{\mu}_k\leftarrow\mu_k$, $\bar{\beta}_k\leftarrow\beta_k$, $k \leftarrow k+1$ and go to \textbf{Step 1}.
\end{itemize}
\textbf{end while} \vspace{1mm} \\
\textbf{Output}: $\bm{x}^k$ \vspace{1mm}
\end{algorithm}

\begin{remark}[\textbf{Comments on special cases of PGels}]\label{eqalgo}
In Algorithm \ref{alg_PGels}, if $\delta=0$, then we have $\beta_k^0 = 0$ and hence $\bm{y}^k = \bm{x}^k$ for all $k\geq0$. In this case, our line search criterion \eqref{lscond} reduces to
\begin{eqnarray*}
F(\bm{u}) - \max\limits_{[k-N]_{+}\leq i\leq k}F(\bm{x}^{i}) \leq -\frac{c}{2}\|\bm{u}-\bm{x}^{k}\|^2.
\end{eqnarray*}
Thus, our PGels readily reduces to the NPG for solving problem \eqref{model} (see, for example, \cite{clp2016penalty,gzlhy2013a,wnf2009sparse}). On the other hand, for any $0<\delta<1$, if
\begin{eqnarray*}
\mu^0_k=\mu_{\max}\geq\frac{L_f+2c}{1-\delta} \quad \mathrm{and} \quad \beta^0_k \leq \sqrt{\frac{\delta(\mu_{\max}-L_f)\mu_{\max}}{4(\mu_{\max}+L_f)^2}} \quad \mathrm{for}~\mathrm{all}~k\geq0,
\end{eqnarray*}
then it follows from Lemma \ref{suffdes} that the line search criterion \eqref{lscond} holds trivially for all $k\geq0$. Thus, in this case, we do not need to perform the line search loop and hence our PGels reduces to a PGe for solving problem \eqref{model}. Finally, if $\delta=0$ and $\mu^0_k=\mu_{\max}\geq L_f + 2c$, it is clear that our PGels reduces to the PG method for solving problem \eqref{model}.
\end{remark}

\begin{remark}[\textbf{Comments on extrapolation parameters in PGels}]
Unlike most of the existing PGe and its variants (mentioned in Section \ref{secintro}) that should choose the extrapolation parameters under certain schemes or conditions, our PGels can choose any $\beta_k^0\in[0,\,\delta\beta_{\max}]$ as an initial guess at each iteration and then updates $\beta_k$ and $\mu_k$ adaptively at the same time if the line search criterion is not satisfied. This strategy actually allows more flexibilities in the choices of the extrapolation parameters, and works well as observed from our computational results in Section \ref{secnum}.
\end{remark}


In the following, we will study the convergence properties of the PGels. Before proceeding, we present the first-order optimality condition for the subproblem \eqref{subpro} in step \textbf{(1b)} of Algorithm \ref{alg_PGels} as follows:
\begin{eqnarray}\label{suboptcond}
0 \in \nabla f(\bm{y}^k) + \mu_k(\bm{u} - \bm{y}^k) + \partial P(\bm{u}).
\end{eqnarray}
We now start our convergence analysis by proving the following lemma, which characterizes the descent property of our potential function.

\begin{lemma}[\textbf{Sufficient descent of $H_{\delta}$}]\label{suffdes}
Suppose that Assumption \ref{assumfun} holds and $\delta \in [0, \,1)$ is a nonnegative constant. Let $\{\bm{x}^k\}$ and $\{\bar{\mu}_k\}$ be the sequences generated by Algorithm \ref{alg_PGels}, and let $\bm{u}$ be the candidate generated by step (1b) at the $k$-th iteration. For any $k \geq 0$, if
\begin{eqnarray*}
\mu_k > L_f \quad \mathrm{and} \quad \beta_k \leq \sqrt{\frac{\delta(\mu_k-L_f)\bar{\mu}_{k-1}}{4(\mu_k+L_f)^2}},
\end{eqnarray*}
then we have
\begin{eqnarray}\label{succhange}
H_{\delta}(\bm{u}, \bm{x}^k, \mu_{k}) - H_{\delta}(\bm{x}^{k}, \bm{x}^{k-1}, \bar{\mu}_{k-1}) \leq - \frac{(1-\delta)\mu_k-L_f}{4}\|\bm{u} - \bm{x}^k\|^2,
\end{eqnarray}
where $H_{\delta}$ is the potential function defined in \eqref{defpofun}.
\end{lemma}
\beginproof
First, from \eqref{subpro}, we have
\begin{eqnarray*}
\langle \nabla f(\bm{y}^k), \,\bm{u} - \bm{y}^k\rangle + \frac{\mu_k}{2}\|\bm{u} - \bm{y}^k\|^2 + P(\bm{u}) \leq
\langle \nabla f(\bm{y}^k), \,\bm{x}^k - \bm{y}^k\rangle + \frac{\mu_k}{2}\|\bm{x}^k - \bm{y}^k\|^2 + P(\bm{x}^k),
\end{eqnarray*}
which implies that
\begin{eqnarray}\label{ineq1}
\begin{aligned}
P(\bm{u})
&\leq P(\bm{x}^k) + \langle \nabla f(\bm{y}^k), \,\bm{x}^k - \bm{u}\rangle + \frac{\mu_k}{2}\|\bm{x}^k - \bm{y}^k\|^2 - \frac{\mu_k}{2}\|\bm{u} - \bm{y}^k\|^2 \\
&= P(\bm{x}^k) + \langle \nabla f(\bm{y}^k), \,\bm{x}^k - \bm{u}\rangle + \frac{\mu_k}{2}\|\bm{x}^k - \bm{y}^k\|^2 - \frac{\mu_k}{2}\|(\bm{u} - \bm{x}^k) + (\bm{x}^k - \bm{y}^k)\|^2 \\
&= P(\bm{x}^k) + \langle \nabla f(\bm{y}^k), \,\bm{x}^k - \bm{u}\rangle - \frac{\mu_k}{2}\|\bm{u} - \bm{x}^k\|^2 + \mu_k \langle \bm{x}^k - \bm{u}, \,\bm{x}^k - \bm{y}^k \rangle
\end{aligned}
\end{eqnarray}
On the other hand, using the fact that $\nabla f$ is Lipschitz continuous with a Lipschitz constant $L_f$ (Assumption \ref{assumfun}(i)), we see from \cite[Lemma 1.2.3]{n2013introductory} that
\begin{eqnarray}\label{ineq2}
f(\bm{u}) \leq f(\bm{x}^k) + \langle \nabla f(\bm{x}^k), \,\bm{u} - \bm{x}^k\rangle + \frac{L_f}{2} \|\bm{u} - \bm{x}^k\|^2.
\end{eqnarray}
Summing \eqref{ineq1} and \eqref{ineq2}, we obtain that
\begin{eqnarray*}
\begin{aligned}
&\quad f(\bm{u}) + P(\bm{u}) - f(\bm{x}^k) - P(\bm{x}^k)\\
&\leq  - \frac{\mu_k-L_f}{2}\|\bm{u} - \bm{x}^k\|^2 + \mu_k \langle \bm{x}^k - \bm{u}, \,\bm{x}^k - \bm{y}^k \rangle
+ \langle \nabla f(\bm{x}^k) - \nabla f(\bm{y}^k), \,\bm{u} - \bm{x}^k\rangle \\
&\leq  - \frac{\mu_k-L_f}{2}\|\bm{u} - \bm{x}^k\|^2 + \left(\mu_k\|\bm{x}^k - \bm{y}^k\|+\|\nabla f(\bm{x}^k) - \nabla f(\bm{y}^k)\|\right)\|\bm{u} - \bm{x}^k\| \\
&\leq  - \frac{\mu_k-L_f}{2}\|\bm{u} - \bm{x}^k\|^2 + (\mu_k+L_f)\|\bm{x}^k - \bm{y}^k\|\|\bm{u} - \bm{x}^k\| \\
&\leq  - \frac{\mu_k-L_f}{2}\|\bm{u} - \bm{x}^k\|^2 + \frac{\mu_k-L_f}{4}\|\bm{u} - \bm{x}^k\|^2 + \frac{(\mu_k+L_f)^2}{\mu_k-L_f}\|\bm{x}^k - \bm{y}^k\|^2 \\
&= - \frac{\mu_k-L_f}{4}\|\bm{u} - \bm{x}^k\|^2 + \frac{(\mu_k+L_f)^2}{\mu_k-L_f}\beta^2_k\|\bm{x}^k - \bm{x}^{k-1}\|^2  \\
&\leq  - \frac{\mu_k-L_f}{4}\|\bm{u} - \bm{x}^k\|^2 + \frac{\delta\bar{\mu}_{k-1}}{4}\|\bm{x}^k - \bm{x}^{k-1}\|^2  \\
&= - \frac{(1-\delta)\mu_k-L_f}{4}\|\bm{u} - \bm{x}^k\|^2 - \frac{\delta\mu_{k}}{4}\|\bm{u} - \bm{x}^k\|^2 + \frac{\delta\bar{\mu}_{k-1}}{4}\|\bm{x}^k - \bm{x}^{k-1}\|^2,
\end{aligned}
\end{eqnarray*}
where the second inequality follows from Cauchy--Schwarz inequality; the third inequality follows from Lipschitz continuity of $\nabla f$; the fourth inequality follows from the relation $ab \leq \frac{a^2}{4s} + s b^2$ with $a=\|\bm{u} - \bm{x}^k\|$, $b=(\mu_k+L_f)\|\bm{x}^k - \bm{y}^k\|$ and $s = \frac{1}{\mu_k-L_f} > 0$; the first equality follows from \eqref{extrastep}; the last inequality follows from $\beta_k \leq \sqrt{\frac{\delta(\mu_k-L_f)\bar{\mu}_{k-1}}{4(\mu_k+L_f)^2}}$. Then, rearranging terms in above relation and recalling the definition of $H_{\delta}$ in \eqref{defpofun}, we obtain \eqref{succhange}.
\endproof

\begin{remark}[\textbf{Comments on Lemma \ref{suffdes}}]
Note that the descent property in Lemma \ref{suffdes} is established for $H_{\delta}$ without requiring $f$ or $P$ to be convex or difference-of-convex function. In fact, with additional assumptions (e.g., convexity) on $f$ or $P$, one can establish a similar descent property for some other constructed potential function $\widetilde{H}$; see, for example, \cite[Lemma 3.1]{wcp2017linear}. Then, one can perform the line search criterion \eqref{lscond} with $\widetilde{H}$ in place of $H_{\delta}$ in Algorithm \ref{alg_PGels} and the convergence analysis can follow in a similar way as presented in this paper. Thus, one can choose suitable potential function in the PGels to fit different scenarios. In this paper, we only focus on $H_{\delta}$ under Assumption \ref{assumfun} to present our main idea.
\end{remark}

From Lemma \ref{suffdes}, we see that the sufficient descent of $H_{\delta}$ can be guaranteed as long as $\mu_k$ is sufficiently large and $\beta_k$ is sufficiently small for each $k\geq0$. Thus, resorting to this lemma, we can show in the following proposition that the line search criterion \eqref{lscond} in Algorithm \ref{alg_PGels} is well-defined.

\begin{proposition}[\textbf{Well-definedness of the line search criterion}]\label{guals}
Suppose that Assumption \ref{assumfun} holds and $\delta \in [0, \,1)$ is a nonnegative constant. Let $\{\bm{x}^k\}$ and $\{\bar{\mu}_k\}$ be the sequences generated by Algorithm \ref{alg_PGels}. Then, for each $k \geq 0$, the line search criterion \eqref{lscond} is satisfied after finitely many inner iterations.
\end{proposition}
\beginproof
We prove this proposition by contradiction. Assume that there exists a $k \geq 0$ such that the line search criterion \eqref{lscond} cannot be satisfied after finitely many inner iterations. Since $\mu_k \leq \mu_{\max}$ due to (1d) in Algorithm \ref{alg_PGels}, then $\mu_k = \mu_{\max}$ must be satisfied after finitely many inner iterations. Let $n_k$ denote the number of inner iterations when $\mu_k = \mu_{\max}$ is satisfied for the \textit{first} time. If $\mu^0_k = \mu_{\max}$, then $n_k=1$; otherwise, we have
\begin{eqnarray*}
\mu_{\min} \tau^{n_k-1} \leq \mu^0_k \tau^{n_k-1} < \mu_{\max},
\end{eqnarray*}
which implies that
\begin{eqnarray*}
n_k \leq \left\lfloor\frac{\log(\mu_{\max})-\log(\mu_{\min})}{\log\tau} + 1 \right\rfloor.
\end{eqnarray*}
Now, we let $\bar{\beta}_k:=\sqrt{\frac{\delta(\mu_{\max}-L_f)\bar{\mu}_{k-1}}{4(\mu_{\max}+L_f)^2}}$ for simplicity. Then, from (1d) in Algorithm \ref{alg_PGels}, we see that $\beta_k$ is decreasing in the inner loop and hence $\beta_k \leq \bar{\beta}_k$ must be satisfied after finitely many inner iterations. Similarly, let $\hat{n}_k$ denote the number of inner iterations when $\beta_k \leq \bar{\beta}_k$ is satisfied for the \textit{first} time. Note that if $\delta=0$, we have $\beta_k^0=\bar{\beta}_k=0$ and hence $\hat{n}_k=1$. For $0<\delta<1$, if $\beta^0_k \leq \bar{\beta}_k$, then $\hat{n}_k \leq n_k$; otherwise, we have
\begin{eqnarray*}
\bar{\beta}_k < \beta^0_k \eta^{\hat{n}_k-1} \leq \delta\beta_{\max} \eta^{\hat{n}_k-1},
\end{eqnarray*}
which implies that
\begin{eqnarray*}
\hat{n}_k \leq \left\lfloor\frac{\log(\delta\beta_{\max})-\log(\bar{\beta}_k)}{-\log\eta} + 1 \right\rfloor.
\end{eqnarray*}
Thus, after at most $\max\{n_k, \,\hat{n}_k\}+1$ inner iterations, we must have $\mu_k \equiv \mu_{\max}$ and $\beta_k \leq \bar{\beta}_k$. Since $c>0$ and $0\leq\delta<1$, one can see that $\mu_{\max}\geq\frac{L_f+2c}{1-\delta} > L_f$ and $(1-\delta)\mu_{\max}-L_f \geq 2c$. Then, using these facts and Lemma \ref{suffdes}, we have
\begin{eqnarray*}
H_{\delta}(\bm{u}, \bm{x}^k, \mu_{\max}) - H_{\delta}(\bm{x}^{k}, \bm{x}^{k-1}, \bar{\mu}_{k-1}) \leq - \frac{(1-\delta)\mu_{\max}-L_f}{4}\|\bm{u} - \bm{x}^k\|^2 \leq - \frac{c}{2}\|\bm{u} - \bm{x}^k\|^2,
\end{eqnarray*}
which, together with
\begin{eqnarray*}
H_{\delta}(\bm{x}^{k}, \bm{x}^{k-1}, \bar{\mu}_{k-1}) \leq \max\limits_{[k-N]_{+}\leq i\leq k}H_{\delta}(\bm{x}^{i}, \bm{x}^{i-1}, \bar{\mu}_{i-1}),
\end{eqnarray*}
implies that \eqref{lscond} must be satisfied after at most $\max\{n_k, \,\hat{n}_k\}+1$ inner iterations. This leads to a contradiction.
\endproof

We are now ready to show our first convergence result in the following theorem that characterizes a cluster point of the sequence generated by the PGels. Our proof is inspired by that of \cite[Lemma 4]{wnf2009sparse}. However, the arguments involved relies on our potential function \eqref{defpofun} that contains multiple blocks of variables. This makes the analysis more intricate. For notational simplicity, from now on, let
\begin{eqnarray}\label{defsk}
\ell(k) \in \mathrm{Arg}\max\limits_{i}\{ H_{\delta}(\bm{x}^{i}, \bm{x}^{i-1}, \bar{\mu}_{i-1}) : i = [k-N]_+, \cdots, k \}.
\end{eqnarray}

\begin{theorem}\label{subconvergence}
Suppose that Assumption \ref{assumfun} holds and $\delta \in [0, \,1)$ is a nonnegative constant. Let $\{\bm{x}^k\}$ and $\{\bar{\mu}_k\}$ be the sequences generated by Algorithm \ref{alg_PGels}. Then,
\begin{itemize}[leftmargin=8mm]
\item[{\rm (i)}] \textbf{(boundedness of sequence)} the sequence $\{\bm{x}^k\}$ is bounded;

\item[{\rm (ii)}] \textbf{(non-increase of subsequence of $H_{\delta}$)} the sequence $\{H_{\delta}(\bm{x}^{\ell(k)}, \bm{x}^{\ell(k)-1}, \bar{\mu}_{\ell(k)-1})\}$ is non-increasing;

\item[{\rm (iii)}] \textbf{(existence of limit)} $\zeta:=\lim\limits_{k\rightarrow\infty} H_{\delta}(\bm{x}^{\ell(k)}, \bm{x}^{\ell(k)-1}, \bar{\mu}_{\ell(k)-1})$ exists;

\item[{\rm (iv)}] \textbf{(diminishing successive changes)} $\lim\limits_{k\rightarrow \infty} \|\bm{x}^{k+1} - \bm{x}^k\| = 0$;

\item[{\rm (v)}] \textbf{(global subsequential convergence)} any cluster point $\bm{x}^*$ of $\{\bm{x}^k\}$ is a stationary point of problem \eqref{model}.
\end{itemize}
\end{theorem}
\beginproof
\textit{Statement (i)}. We first prove by induction that
\begin{eqnarray}\label{bdbyin}
H_{\delta}(\bm{x}^{k}, \bm{x}^{k-1}, \bar{\mu}_{k-1}) \leq F(\bm{x}^0)
\end{eqnarray}
for all $k \geq 1$. Indeed, for $k=1$, it follows from Proposition \ref{guals} that
\begin{eqnarray*}
H_{\delta}(\bm{x}^1, \bm{x}^0, \bar{\mu}_{0}) - H_{\delta}(\bm{x}^{0}, \bm{x}^{-1}, \bar{\mu}_{-1}) \leq -\frac{c}{2} \|\bm{x}^1-\bm{x}^0\|^2 \leq 0
\end{eqnarray*}
is satisfied after finitely many inner iterations. This, together with $\bm{x}^{-1}=\bm{x}^0$, implies that
\begin{eqnarray*}
H_{\delta}(\bm{x}^1, \bm{x}^0, \bar{\mu}_{0}) \leq H_{\delta}(\bm{x}^{0}, \bm{x}^{-1}, \bar{\mu}_{-1}) = F(\bm{x}^0).
\end{eqnarray*}
Hence, \eqref{bdbyin} holds for $k=1$. We now suppose that \eqref{bdbyin} holds for all $k \leq K$ for some integer $K \geq 1$. Next, we show that \eqref{bdbyin} also holds for $k=K+1$. Indeed, for $k=K+1$, we have
\begin{eqnarray*}
\begin{aligned}
H_{\delta}(\bm{x}^{K+1}, \bm{x}^{K}, \bar{\mu}_{K}) - F(\bm{x}^0)
&\leq H_{\delta}(\bm{x}^{K+1}, \bm{x}^{K}, \bar{\mu}_{K}) - \max\limits_{[K-N]_{+}\leq i\leq K} H_{\delta}(\bm{x}^{i}, \bm{x}^{i-1}, \bar{\mu}_{i-1})   \\
&\leq -\frac{c}{2}\|\bm{x}^{K+1}-\bm{x}^{K}\|^2 \leq 0,
\end{aligned}
\end{eqnarray*}
where the first inequality follows from the induction hypothesis and the second inequality follows from \eqref{lscond}. Hence, \eqref{bdbyin} holds for $k=K+1$. This completes the induction. Then, from \eqref{bdbyin}, we have that for any $k \geq 1$,
\begin{eqnarray*}
F(\bm{x}^0) \geq H_{\delta}(\bm{x}^{k}, \bm{x}^{k-1}, \bar{\mu}_{k-1}) \geq F(\bm{x}^k),
\end{eqnarray*}
which, together with Assumption \ref{assumfun}(iii), implies that $\{\bm{x}^k\}$ is bounded. This proves statement (i).

\textit{Statement (ii)}.
Recall the definition of $\ell(k)$ in \eqref{defsk} and let $\Delta_{\bm{x}^{k}}:=\bm{x}^{k+1}-\bm{x}^{k}$ for simplicity. Then, from the line search criterion \eqref{lscond}, we have
\begin{eqnarray}\label{lscond2}
H_{\delta}(\bm{x}^{k+1}, \bm{x}^{k}, \bar{\mu}_{k}) - H_{\delta}(\bm{x}^{\ell(k)}, \bm{x}^{\ell(k)-1}, \bar{\mu}_{\ell(k)-1}) \leq -\frac{c}{2}\|\Delta_{\bm{x}^{k}}\|^2 \leq 0.
\end{eqnarray}
Observe that
\begin{eqnarray*}
\begin{aligned}
&\quad H_{\delta}(\bm{x}^{\ell(k+1)}, \bm{x}^{\ell(k+1)-1}, \bar{\mu}_{\ell(k+1)-1})   \\
&= \max\limits_{[k+1-N]_{+}\leq i\leq k+1} H_{\delta}(\bm{x}^{i}, x^{i-1}, \bar{\mu}_{i-1})  \\
&= \max\left\{H_{\delta}(\bm{x}^{k+1}, \bm{x}^{k}, \bar{\mu}_{k}), \,\max\limits_{[k+1-N]_{+}\leq i\leq k} H_{\delta}(\bm{x}^{i}, \bm{x}^{i-1}, \bar{\mu}_{i-1})\right\} \\
&\leq \max\left\{H_{\delta}(\bm{x}^{\ell(k)}, \bm{x}^{\ell(k)-1}, \bar{\mu}_{\ell(k)-1}), \,\max\limits_{[k+1-N]_{+}\leq i\leq k} H_{\delta}(\bm{x}^{i}, \bm{x}^{i-1}, \bar{\mu}_{i-1})\right\}    \\
&\leq \max\left\{H_{\delta}(\bm{x}^{\ell(k)}, \bm{x}^{\ell(k)-1}, \bar{\mu}_{\ell(k)-1}), \,\max\limits_{[k-N]_{+}\leq i\leq k} H_{\delta}(\bm{x}^{i}, \bm{x}^{i-1}, \bar{\mu}_{i-1})\right\}    \\
&= \max\left\{H_{\delta}(\bm{x}^{\ell(k)}, \bm{x}^{\ell(k)-1}, \bar{\mu}_{\ell(k)-1}), \,H_{\delta}(\bm{x}^{\ell(k)}, \bm{x}^{\ell(k)-1}, \bar{\mu}_{\ell(k)-1})\right\} \\
&= H_{\delta}(\bm{x}^{\ell(k)}, \bm{x}^{\ell(k)-1}, \bar{\mu}_{\ell(k)-1}),
\end{aligned}
\end{eqnarray*}
where the first inequality follows from \eqref{lscond2} and the second last equality follows from \eqref{defsk}. This proves statement (ii).

\textit{Statement (iii)}.
It follows from Assumption \ref{assumfun}(iii) and the definition of $H_{\delta}$ in \eqref{defpofun} that $H_{\delta}(\bm{x}^{\ell(k)}$, $\bm{x}^{\ell(k)-1}, \bar{\mu}_{\ell(k)-1})$ is bounded below. This together with statement (ii) proves that there exists a number $\zeta$ such that
\begin{eqnarray}\label{limsk}
\lim\limits_{k\rightarrow\infty} H_{\delta}(\bm{x}^{\ell(k)}, \bm{x}^{\ell(k)-1}, \bar{\mu}_{\ell(k)-1}) = \zeta.
\end{eqnarray}

\textit{Statement (iv)}.
We next prove statement (iv). To this end, we first show by induction that for all $j\geq1$, it holds that
\begin{numcases}{}
\lim\limits_{k\rightarrow\infty} \Delta_{\bm{x}^{\ell(k)-j}} = 0,  \label{keylim1} \\
\lim\limits_{k\rightarrow\infty} F(\bm{x}^{\ell(k)-j}) = \zeta \label{keylim2}.
\end{numcases}
We start by proving \eqref{keylim1} and \eqref{keylim2} for $j=1$. Applying \eqref{lscond2} with $k$ replaced by $\ell(k)-1$, we obtain
\begin{eqnarray*}
H_{\delta}(\bm{x}^{\ell(k)}, \bm{x}^{\ell(k)-1}, \bar{\mu}_{\ell(k)-1}) - H_{\delta}(\bm{x}^{\ell(\ell(k)-1)}, \bm{x}^{\ell(\ell(k)-1)-1}, \bar{\mu}_{\ell(\ell(k)-1)-1}) \leq -\frac{c}{2}\|\Delta_{\bm{x}^{\ell(k)-1}}\|^2,
\end{eqnarray*}
which, together with \eqref{limsk}, implies that
\begin{eqnarray}\label{limsk1}
\lim\limits_{k\rightarrow\infty} \Delta_{\bm{x}^{\ell(k)-1}} = 0.
\end{eqnarray}
Then, from \eqref{limsk} and \eqref{limsk1}, we have
\begin{eqnarray*}
\begin{aligned}
\zeta
&=\lim\limits_{k\rightarrow\infty} H_{\delta}(\bm{x}^{\ell(k)}, \bm{x}^{\ell(k)-1}, \bar{\mu}_{\ell(k)-1})  \\
&=\lim\limits_{k\rightarrow\infty} F(\bm{x}^{\ell(k)-1}+\Delta_{\bm{x}^{\ell(k)-1}}) + \frac{\delta\bar{\mu}_{\ell(k)-1}}{4}\|\Delta_{\bm{x}^{\ell(k)-1}}\|^2 \\
&=\lim\limits_{k\rightarrow\infty} F(\bm{x}^{\ell(k)-1}),
\end{aligned}
\end{eqnarray*}
where the second equality follows from the definition of $H_{\delta}$ in \eqref{defpofun} and the last equality follows because $\{\bm{x}^k\}$ is bounded (see statement (i)), $\{\bar{\mu}_k\}$ is bounded (since $\mu_{\min}\leq\bar{\mu}_k\leq\mu_{\max}$) and $F$ is uniformly continuous on any compact subset of $\mathrm{dom}\,F$ under Assumption \ref{assumfun}(i) and (ii). Thus, \eqref{keylim1} and \eqref{keylim2} hold for $j=1$.

We next suppose that \eqref{keylim1} and \eqref{keylim2} hold for $j=J$ for some $J \geq 1$. It remains to show that they also hold for $j = J+1$. Indeed, from \eqref{lscond2} with $k$ replaced by $\ell(k)-J-1$ (here, without loss of generality, we assume that $k$ is large enough such that $\ell(k)-J-1$ is nonnegative), we have
\begin{eqnarray*}
\begin{aligned}
&~~H_{\delta}(\bm{x}^{\ell(k)-J}, \bm{x}^{\ell(k)-J-1}, \bar{\mu}_{\ell(k)-J-1}) - H_{\delta}(\bm{x}^{\ell(\ell(k)-J-1)}, \bm{x}^{\ell(\ell(k)-J-1)-1}, \bar{\mu}_{\ell(\ell(k)-J-1)-1}) \\
&\leq -\frac{c}{2}\|\Delta_{\bm{x}^{\ell(k)-J-1}}\|^2,
\end{aligned}
\end{eqnarray*}
which, together with the definition of $H_{\delta}$ in \eqref{defpofun}, implies that
\begin{eqnarray*}
\left(\frac{c}{2}+\frac{\delta\bar{\mu}_{\ell(k)-J-1}}{4}\right)\|\Delta_{\bm{x}^{\ell(k)-J-1}}\|^2 \leq H_{\delta}(\bm{x}^{\ell(\ell(k)-J-1)}, \bm{x}^{\ell(\ell(k)-J-1)-1}, \bar{\mu}_{\ell(\ell(k)-J-1)-1})- F(\bm{x}^{\ell(k)-J}).
\end{eqnarray*}
This together with \eqref{limsk} and the induction hypothesis implies that
\begin{eqnarray*}
\lim\limits_{k\rightarrow\infty} \Delta_{\bm{x}^{\ell(k)-(J+1)}} = 0.
\end{eqnarray*}
Thus, \eqref{keylim1} holds for $j=J+1$. From this, we further have
\begin{eqnarray*}
\lim\limits_{k\rightarrow\infty} F(\bm{x}^{\ell(k)-(J+1)})
= \lim\limits_{k\rightarrow\infty} F(\bm{x}^{\ell(k)-J} - \Delta_{\bm{x}^{\ell(k)-(J+1)}})
= \lim\limits_{k\rightarrow\infty} F(\bm{x}^{\ell(k)-J}) = \zeta,
\end{eqnarray*}
where the second equality follows because $\{\bm{x}^k\}$ is bounded (see statement (i)), $\{\bar{\mu}_k\}$ is bounded (since $\mu_{\min}\leq\bar{\mu}_k\leq\mu_{\max}$) and $F$ is uniformly continuous on any compact subset of $\mathrm{dom}\,F$ under Assumption \ref{assumfun}(i) and (ii). Hence, \eqref{keylim2} also holds for $j=J+1$. This completes the induction.

We are now ready to prove the main result in this statement. Indeed, recalling the definition of $\ell(k)$ in \eqref{defsk}, we see that $k-N\leq \ell(k)\leq k$ (without loss of generality, we assume that $k$ is large enough such that $k \geq N$). Thus, for any $k$, we must have $k-N-1=\ell(k)-j_k$ for some $j_k \in [1,\,N+1]$. Then, we have
\begin{eqnarray*}
\|\Delta_{\bm{x}^{k-N-1}}\|=\|\Delta_{\bm{x}^{\ell(k)-j_k}}\| \leq \max\limits_{1 \leq j \leq N+1} \|\Delta_{\bm{x}^{\ell(k)-j}}\|.
\end{eqnarray*}
This together with \eqref{keylim1} implies that
\begin{eqnarray*}
\lim\limits_{k\rightarrow\infty} \Delta_{\bm{x}^{k}} = \lim\limits_{k\rightarrow\infty} \Delta_{\bm{x}^{k-N-1}} = 0,
\end{eqnarray*}
This proves statement (iv).

\textit{Statement (v)}. First, since $\{\bm{x}^k\}$ is bounded (see statement (i)), there exists at least one cluster point. Suppose that $\bm{x}^*$ is a cluster point of $\{\bm{x}^k\}$ and let $\{\bm{x}^{k_i}\}$ be a convergent subsequence such that $\lim\limits_{i\rightarrow\infty} \bm{x}^{k_i} = \bm{x}^*$. Then, we see from statement (iv) and the boundedness of $\beta_k$ for any $k$ (since $0\leq\beta_k\leq\beta^0_k\leq\delta\beta_{\max}$) that
\begin{eqnarray}\label{limy}
\lim\limits_{i\rightarrow\infty} \bm{y}^{k_i} = \lim\limits_{i\rightarrow\infty} \bm{x}^{k_i} + \beta_{k_i}(\bm{x}^{k_i} - \bm{x}^{k_i-1})
= \lim\limits_{i\rightarrow\infty} \bm{x}^{k_i} = \bm{x}^*.
\end{eqnarray}
Thus, passing to the limit along $\{(\bm{x}^{k_i}, \,\bm{y}^{k_i})\}$ in \eqref{suboptcond} with $\bm{x}^{k_i+1}$ in place of $\bm{u}$ and $\bar{\mu}_{k_i}$ in place of $\mu_{k}$, and invoking Assumption \ref{assumfun}(ii), statement (iv), the boundedness of $\{\bar{\mu}_k\}$ (since $\mu_{\min}\leq\bar{\mu}_k\leq\mu_{\max}$), \eqref{robust} and \eqref{limy}, we obtain
\begin{eqnarray*}
0 \in \nabla f(\bm{x}^*) + \partial P(\bm{x}^*),
\end{eqnarray*}
which implies that $\bm{x}^*$ is a stationary point of \eqref{model}. This proves statement (v).
\endproof

Based on Theorem \ref{subconvergence}, we can further characterize the sequence of objective values along $\{\bm{x}^k\}$ in the following proposition. This proposition will be useful in the analysis of the local convergence rate in the next section.

\begin{proposition}\label{subfunseq}
Suppose that Assumption \ref{assumfun} holds and $\delta \in [0, \,1)$ is a nonnegative constant. Let $\{\bm{x}^k\}$ be a sequence generated by Algorithm \ref{alg_PGels} and $\ell(k)$ be the index defined in \eqref{defsk} for each $k$. Then,
\begin{itemize}
\item[{\rm (i)}] $\lim\limits_{k\rightarrow\infty} F(\bm{x}^{\ell(k)}) = \zeta$, where $\zeta$ is given in Theorem \ref{subconvergence}(iii);

\item[{\rm (ii)}] $F \equiv \zeta$ on $\Omega$, where $\Omega$ is the set of cluster points of the subsequence $\{\bm{x}^{\ell(k)}\}$.
\end{itemize}
\end{proposition}
\beginproof
Statement (i) follows immediately from Theorem \ref{subconvergence}(iii), $\|\bm{x}^{k+1}-\bm{x}^k\|\rightarrow0$ (see Theorem \ref{subconvergence}(iv)), the boundedness of $\{\bar{\mu}_k\}$ (since $\mu_{\min}\leq\bar{\mu}_k\leq\mu_{\max}$) and the definition of $H_{\delta}$ in \eqref{defpofun}.

We now prove statement (ii). First, since $\{\bm{x}^{\ell(k)}\}$ is a subsequence of $\{\bm{x}^k\}$, it follows from Theorem \ref{subconvergence}(i) and (v) that $\emptyset \neq \Omega \subseteq\mathcal{X}$, where $\mathcal{X}$ is the set of all stationary points of \eqref{model}. Moreover, we recall from Assumption \ref{assumfun}(i) and (ii) that $f$ is continuously differentiable and $P$ is continuous on its domain. These facts prove statement (ii).
\endproof

\section{Local convergence rate of two special cases of PGels}\label{secconrate}

In this section, under the additional assumption that the objective $F$ in problem \eqref{model} is a KL function with an exponent $\theta$, we further study the local convergence rate of two special cases of the PGels in terms of objective values. The first case is for the PGels with $\delta=0$, namely, the NPG and the other is for the PGels with $N=0$, where $N$ is the line search gap used in \eqref{lscond}.

\subsection{Local convergence rate of NPG}

In this subsection, we discuss the local convergence rate of the PGels with $\delta=0$, namely, the NPG (see Remark \ref{eqalgo}). In this case, we have $\bm{y}^k \equiv \bm{x}^k$ and $F(\bm{x}^{k}) \equiv H_{\delta}(\bm{x}^{k}, \bm{x}^{k-1}, \bar{\mu}_{k-1})$ for all $k\geq0$. The main results are presented in the following theorem. This kind of results on local convergence rate have been well studied for many existing algorithms; see, for example, \cite{ab2009on,fgp2015splitting,wcp2017a}. However, the analysis there heavily relies on the \textit{monotonicity} of the objective or certain potential function along the sequence generated and hence cannot be applied for the NPG. More intricate analysis is needed for handling the non-monotone line search. To the best of our knowledge, this is the first local convergence rate analysis of the NPG for solving \eqref{model}. This analysis framework may also be applied for some other non-monotone algorithms, e.g., a non-monotone alternating updating method for a class of matrix factorization problems \cite{ypc2017a} (see also \cite[Theorem 4.4]{yang2017first}).

\begin{theorem}\label{conratenpg}
Suppose that Assumption \ref{assumfun} holds and the objective $F$ in problem \eqref{model} is a KL function with an exponent $\theta$. Let $\{\bm{x}^k\}$ and $\{\bar{\mu}_k\}$ be the sequences generated by Algorithm \ref{alg_PGels} with $\delta=0$, and let $\zeta$ be given in Theorem \ref{subconvergence}(iii). Then, the following statements hold.
\begin{itemize}
\item[{\rm (i)}] If $\theta = 0$, $F(\bm{x}^k) \leq \zeta$ for all large $k$;
\item[{\rm (ii)}] If $\theta \in \left(0, \frac{1}{2}\right]$, there exist $\rho\in(0,1)$ and $c_1>0$ such that $F(\bm{x}^k) - \zeta \leq c_1\rho^k$ for all large $k$;
\item[{\rm (iii)}] If $\theta \in \left(\frac{1}{2}, 1\right)$, there exists $c_2>0$ such that $F(\bm{x}^k) - \zeta \leq c_2\,k^{-\frac{1}{2\theta-1}}$ for all large $k$.
\end{itemize}
\end{theorem}
\beginproof
We start by defining an index sequence $\{\xi(t)\}_{t=0}^{\infty}$ as follows:
\begin{eqnarray*}
\xi(t)=\ell((N+1)t), \quad t = 0,\,1,\,2,\,\cdots,
\end{eqnarray*}
where $\ell(k)$ is defined in \eqref{defsk}. It is obvious that $\{\xi(t)\}_{t=0}^{\infty}$ is a subsequence of $\{\ell(k)\}_{k=0}^{\infty}$. Moreover, since $k-N\leq\ell(k)\leq k$ for any $k\geq N$, we have
\begin{eqnarray*}
\xi(t) = \ell((N+1)t) \geq (N+1)t - N = (N+1)(t-1)+1 \geq \ell((N+1)(t-1))+1 > \xi(t-1)
\end{eqnarray*}
for any $t \geq 1$. Therefore, $\{\xi(t)\}_{t=0}^{\infty}$ is increasing. We now recall from Theorem \ref{subconvergence}(ii) and Assumption \ref{assumfun}(iii) that $\{F(\bm{x}^{\ell(k)})\}_{k=0}^{\infty}$ is non-increasing and bounded below. Since $\{F(\bm{x}^{\xi(t)})\}_{t=0}^{\infty}$ is a subsequence of $\{F(\bm{x}^{\ell(k)})\}_{k=0}^{\infty}$, then it follows that $\{F(\bm{x}^{\xi(t)})\}_{t=0}^{\infty}$ is non-increasing and bounded below. Moreover, from Proposition \ref{subfunseq}, we have
\begin{eqnarray*}
\lim\limits_{t\rightarrow\infty} F(\bm{x}^{\xi(t)})=\lim\limits_{k\rightarrow\infty} F(\bm{x}^{\ell(k)})=\zeta.
\end{eqnarray*}
In addition, it follows from \eqref{lscond2} with $k$ replaced by $\xi(t)-1$ that
\begin{eqnarray}\label{suffdesnewnpg}
\begin{aligned}
F(\bm{x}^{\xi(t)})
&\leq F(\bm{x}^{\ell(\xi(t)-1)})-\frac{c}{2}\|\Delta_{\bm{x}^{\xi(t)-1}}\|^2
\leq F(\bm{x}^{\ell((N+1)(t-1))})-\frac{c}{2}\|\Delta_{\bm{x}^{\xi(t)-1}}\|^2 \\
&=F(\bm{x}^{\xi(t-1)})-\frac{c}{2}\|\Delta_{\bm{x}^{\xi(t)-1}}\|^2,
\end{aligned}
\end{eqnarray}
where the second inequality follows because $\{F(\bm{x}^{\ell(k)})\}_{k=0}^{\infty}$ is non-increasing and $\xi(t)-1 = \ell((N+1)t) - 1 \geq (N+1)(t-1)$. We next consider two cases.

\textbf{Case 1.} In this case, we suppose that $F(\bm{x}^{\xi(T)}) = \zeta$ for some $T \geq 0$. Since the sequence $\{F(\bm{x}^{\xi(t)})\}_{t=0}^{\infty}$ is non-increasing, we must have $F(\bm{x}^{\xi(t)}) = \zeta$ for all $t \geq T$. Then, for all $k \in [(N+1)t-N, \,(N+1)t]$ with any $t \geq T$, we have
\begin{eqnarray*}
F(\bm{x}^k) \leq F(\bm{x}^{\ell((N+1)t)}) = F(\bm{x}^{\xi(t)}) = \zeta.
\end{eqnarray*}
Thus, the conclusions of three statements hold.

\textbf{Case 2.} From now on, we consider the case where $F(\bm{x}^{\xi(t)}) > \zeta$ for all $t \geq 0$. From Theorem \ref{subconvergence}(i), we see that $\{\bm{x}^{\xi(t)}\}_{t=0}^{\infty}$ is bounded and hence must have at least one cluster point. Let $\Gamma$ denote the set of cluster points of $\{\bm{x}^{\xi(t)}\}_{t=0}^{\infty}$. Since $\{\bm{x}^{\xi(t)}\}_{t=0}^{\infty}$ is a subsequence of $\{\bm{x}^{\ell(k)}\}_{k=0}^{\infty}$, we have $\Gamma \subseteq \Omega$, where $\Omega$ is the set of cluster points of $\{\bm{x}^{\ell(k)}\}_{k=0}^{\infty}$. Then, it follows from Proposition \ref{subfunseq} that $F \equiv \zeta$ on $\Gamma$. This fact together with our assumption that $F$ is a KL function with an exponent $\theta$ and Proposition \ref{uniKL}, there exist $\varepsilon, \nu>0$ such that
\begin{eqnarray*}
\varphi'\left(F(\bm{x})-\zeta\right)\mathrm{dist}\left(0,
\,\partial F(\bm{x})\right) \geq 1, \quad \mathrm{where}~\varphi(s)=as^{1-\theta}~\mathrm{for}~\mathrm{some}~a>0,
\end{eqnarray*}
for all $\bm{x}$ satisfying $\mathrm{dist}(\bm{x}, \,\Gamma) < \varepsilon$ and $\zeta < F(\bm{x}) < \zeta + \nu$. On the other hand, since $\lim\limits_{t\rightarrow \infty}\mathrm{dist}(\bm{x}^{\xi(t)},\,\Gamma)=0$ (by the definition of $\Gamma$) and $F(\bm{x}^{\xi(t)})\rightarrow\zeta$, then for such $\varepsilon$ and $\nu$, there exists an integer $T_0\geq0$ such that $\mathrm{dist}(\bm{x}^{\xi(t)},\,\Gamma) < \varepsilon$ and $\zeta<F(\bm{x}^{\xi(t)})<\zeta+\nu$ for all $t\geq T_0$. Thus, for $t\geq T_0$, we have
\begin{eqnarray}\label{klineqnpg}
\varphi'\left(F(\bm{x}^{\xi(t)})-\zeta\right)\mathrm{dist}\left(0,~\partial F(\bm{x}^{\xi(t)})\right) \geq 1, \quad \mathrm{where}~\varphi(s)=as^{1-\theta}~\mathrm{for}~\mathrm{some}~a>0.
\end{eqnarray}

Next, looking at the subdifferential $\partial F(\bm{x}^{k})$, we have
\begin{eqnarray*}
\begin{aligned}
\partial F(\bm{x}^{k})
&= \nabla f(\bm{x}^{k}) + \partial P(\bm{x}^{k}) \\
&= \nabla f(\bm{x}^{k-1}) + \partial P(\bm{x}^{k}) + \bar{\mu}_{k-1}(\bm{x}^{k} - \bm{x}^{k-1}) + \nabla f(\bm{x}^{k}) - \nabla f(\bm{x}^{k-1}) -\bar{\mu}_{k-1}(\bm{x}^{k} - \bm{x}^{k-1})   \\
&\ni \nabla f(\bm{x}^{k}) - \nabla f(\bm{x}^{k-1}) - \bar{\mu}_{k-1}(\bm{x}^{k} - \bm{x}^{k-1}),
\end{aligned}
\end{eqnarray*}
where the inclusion follows from the optimality condition for \eqref{subpro} at the $(k\!-\!1)$-st iteration, i.e., $0 \in \nabla f(\bm{x}^{k-1}) + \partial P(\bm{x}^{k}) + \bar{\mu}_{k-1}(\bm{x}^{k} - \bm{x}^{k-1}) $. Using this relation together with the global Lipschitz continuity of $\nabla f$ and the boundedness of $\{\bar{\mu}_k\}$, there exists $d_1>0$ such that
\begin{eqnarray}\label{upbdpatnpg}
\mathrm{dist}\left(0,\,\partial F(\bm{x}^{k})\right) \leq d_1\|\bm{x}^{k}-\bm{x}^{k-1}\|.
\end{eqnarray}

Now, for notational simplicity, let $\Delta_{F}^{\xi(t)}:=F(\bm{x}^{\xi(t)})-\zeta$. Since $\{F(\bm{x}^{\xi(t)})\}_{t=0}^{\infty}$ is non-increasing, we see that $\Delta_{F}^{\xi(t)}$ is non-increasing, $\Delta_{F}^{\xi(t)}>0$ for $t\geq0$ and $\Delta_{F}^{\xi(t)}\rightarrow0$. Then, for all $t \geq T_0$, we have
\begin{eqnarray}\label{keyineqnpg}
\begin{aligned}
1
&\leq \varphi'\left(\Delta_{F}^{\xi(t)}\right)\mathrm{dist}\left( 0,\,\partial F(\bm{x}^{\xi(t)})\right) \\
&\leq \frac{a}{1-\theta}\cdot\left(\Delta_{F}^{\xi(t)}\right)^{-\theta}\cdot d_1\|\bm{x}^{\xi(t)}-\bm{x}^{\xi(t)-1}\| \\
&\leq \frac{ad_1}{1-\theta}\cdot\left(\Delta_{F}^{\xi(t)}\right)^{-\theta}\cdot
\sqrt{\frac{2}{c}\left(F(\bm{x}^{\xi(t-1)})-F(\bm{x}^{\xi(t)})\right)}  \\
&= d_2 \left(\Delta_{F}^{\xi(t)}\right)^{-\theta}\cdot \sqrt{\Delta_{F}^{\xi(t-1)}-\Delta_{F}^{\xi(t)}},
\end{aligned}
\end{eqnarray}
where $d_2:=\frac{\sqrt{2}ad_1}{(1-\theta)\sqrt{c}}>0$, the first inequality follows from \eqref{klineqnpg}, the second inequality follows from \eqref{upbdpatnpg} and the last inequality follows from \eqref{suffdesnewnpg}. Next, we consider the following three cases.

\textbf{(i)} $\theta=0$. In this case, we see from \eqref{keyineqnpg} that $\Delta_{F}^{\xi(t-1)}-\Delta_{F}^{\xi(t)} \geq \frac{1}{d_2^2}$ for all $t \geq T_0$, which contradicts $\Delta_{F}^{\xi(t)}\rightarrow0$. Thus, \textbf{Case 2} cannot happen.

\textbf{(ii)} $0<\theta\leq\frac{1}{2}$. In this case, we have $0<2\theta\leq1$. Since $\Delta_{F}^{\xi(t)}\rightarrow0$, there exists $T_1\geq0$ such that $\Delta_{F}^{\xi(t)} \leq 1$ for all $t \geq \widetilde{T}:=\max\{T_0,\,T_1\}$. Then, for all $t \geq \widetilde{T}$, we see from \eqref{keyineqnpg} that
\begin{eqnarray*}
\Delta_{F}^{\xi(t)}\leq\left(\Delta_{F}^{\xi(t)}\right)^{2\theta} \leq d_2^2\left(\Delta_{F}^{\xi(t-1)}-\Delta_{F}^{\xi(t)}\right),
\end{eqnarray*}
which implies that
\begin{eqnarray*}
\Delta_{F}^{\xi(t)} \leq \gamma \,\Delta_{F}^{\xi(t-1)}
\leq \cdots \leq \gamma^{t-\widetilde{T}+1} \,\Delta_{F}^{\xi(\widetilde{T}-1)},
\end{eqnarray*}
where $\gamma:=\frac{d_2^2}{1+d_2^2} < 1$. Then, for all $k \in [(N+1)t-N, \,(N+1)t]$ with any $t \geq \widetilde{T}$, we have
\begin{eqnarray*}
F(\bm{x}^{k})-\zeta
\leq \Delta_{F}^{\xi(t)}
\leq \gamma^{t-\widetilde{T}+1} \,\Delta_{F}^{\xi(\widetilde{T}-1)}
\leq \gamma^{\frac{k}{N+1}-\widetilde{T}+1} \,\Delta_{F}^{\xi(\widetilde{T}-1)}
= c_1 \rho^k,
\end{eqnarray*}
where $c_1:=\gamma^{-\widetilde{T}+1} \,\Delta_{F}^{\xi(\widetilde{T}-1)}$, $\rho:=\gamma^{\frac{1}{N+1}}<1$ and the last inequality follows from $t \geq \frac{k}{N+1}$. This proves statement (ii).

\textbf{(iii)} $\frac{1}{2}<\theta<1$. We define $g(s):=s^{-2\theta}$ for $s\in(0,\,\infty)$. It is easy to see that $g$ is non-increasing. Then, for any $t \geq T_0$, we further consider the following two cases.
\begin{itemize}
\item If $g(\Delta_{F}^{\xi(t)})\leq 2g(\Delta_{F}^{\xi(t-1)})$, it follows from \eqref{keyineqnpg} that
      \begin{eqnarray*}
      \begin{aligned}
      \frac{1}{d^2_2}
      &\leq (\Delta_{F}^{\xi(t)})^{-2\theta}\cdot(\Delta_{F}^{\xi(t-1)}-\Delta_{F}^{\xi(t)})
       = g(\Delta_{F}^{\xi(t)})\cdot (\Delta_{F}^{\xi(t-1)}-\Delta_{F}^{\xi(t)})
       \leq 2g(\Delta_{F}^{\xi(t-1)})\cdot (\Delta_{F}^{\xi(t-1)}-\Delta_{F}^{\xi(t)})   \\
      &\leq 2\int^{\Delta_{F}^{\xi(t-1)}}_{\Delta_{F}^{\xi(t)}}\!\!\! g(s) \,\mathrm{d}s
       =\frac{2(\Delta_{F}^{\xi(t-1)})^{1-2\theta}-2(\Delta_{F}^{\xi(t)})^{1-2\theta}}{1-2\theta},
      \end{aligned}
      \end{eqnarray*}
      which, together with $1-2\theta<0$, implies that
      \begin{eqnarray}\label{difflbd1npg}
      (\Delta_{F}^{\xi(t)})^{1-2\theta}-(\Delta_{F}^{\xi(t-1)})^{1-2\theta} \geq \frac{2\theta-1}{2d_2^2}.
      \end{eqnarray}

\item If $g(\Delta_{F}^{\xi(t)})\geq 2g(\Delta_{F}^{\xi(t-1)})$, one can check that $(\Delta_{F}^{\xi(t)})^{1-2\theta}\geq2^{\frac{2\theta-1}{2\theta}}(\Delta_{F}^{\xi(t-1)})^{1-2\theta}$. Then, we have
      \begin{eqnarray}\label{difflbd2npg}
      (\Delta_{F}^{\xi(t)})^{1-2\theta}-(\Delta_{F}^{\xi(t-1)})^{1-2\theta}
      \geq \left(2^{\frac{2\theta-1}{2\theta}}-1\right)(\Delta_{F}^{\xi(t-1)})^{1-2\theta}
      \geq \left(2^{\frac{2\theta-1}{2\theta}}-1\right)(\Delta_{F}^{\xi(T_0-1)})^{1-2\theta},
      \end{eqnarray}
      where the last inequality follows from the facts that $\Delta_{F}^{\xi(t)}$ is non-increasing and $1-2\theta<0$.
\end{itemize}

Thus, combining \eqref{difflbd1npg} and \eqref{difflbd2npg}, we obtain
\begin{eqnarray*}
(\Delta_{F}^{\xi(t)})^{1-2\theta}-(\Delta_{F}^{\xi(t-1)})^{1-2\theta}
\geq d_3 := \min\left\{\frac{2\theta-1}{2d_2^2},\,
\left(2^{\frac{2\theta-1}{2\theta}}-1\right)(\Delta_{F}^{\xi(T_0-1)})^{1-2\theta}\right\}.
\end{eqnarray*}
Then, we have
\begin{eqnarray*}
\begin{aligned}
(\Delta_{F}^{\xi(t)})^{1-2\theta}
\geq (\Delta_{F}^{\xi(t)})^{1-2\theta}-(\Delta_{F}^{\xi(T_0)})^{1-2\theta}
=\sum^{t}_{j = T_0+1} \left((\Delta_{F}^{\xi(j)})^{1-2\theta}-(\Delta_{F}^{\xi(j-1)})^{1-2\theta}\right)
\geq (t-T_0)d_3 \geq \frac{d_3}{2}\,t,
\end{aligned}
\end{eqnarray*}
where the last inequality holds for $t \geq 2T_0$. This implies that $\Delta_{F}^{\xi(t)} \leq  d_4\,t^{-\frac{1}{2\theta-1}}$, where $d_4:=\left(\frac{d_3}{2}\right)^{-\frac{1}{2\theta-1}}$. Then, for all $k \in [(N+1)t-N, \,(N+1)t]$ with any $t \geq 2T_0$, we have
\begin{eqnarray*}
F(\bm{x}^{k}) - \zeta \leq \Delta_{F}^{\xi(t)} \leq d_4\,t^{-\frac{1}{2\theta-1}}
\leq c_2\,k^{-\frac{1}{2\theta-1}},
\end{eqnarray*}
where $c_2:=(N+1)^{\frac{1}{2\theta-1}}d_4$ and the last inequality follows from $t \geq \frac{k}{N+1}$. This proves statement (iii).
\endproof

\subsection{Local convergence rate of PGels with $N=0$}

In this subsection, we discuss the local convergence rate of the PGels with $N=0$, i.e., we set the line search gap $N$ to 0 in \eqref{lscond}. We first give the following lemma.

\begin{lemma}\label{pabound}
Suppose that Assumption \ref{assumfun} holds and $\delta \in [0, \,1)$ is a nonnegative constant. Let $\{\bm{x}^k\}$, $\{\bar{\mu}_k\}$ and $\{\bar{\beta}_k\}$ be the sequences generated by Algorithm \ref{alg_PGels}. Then, there exist $\tilde{c}>0$ and $K>0$ such that
\begin{eqnarray}\label{distbd}
\mathrm{dist}\left((0,0,0),\,\partial H_{\delta}(\bm{x}^{k}, \bm{x}^{k-1}, \bar{\mu}_{k-1})\right) \leq \tilde{c}\left(\|\bm{x}^{k}-\bm{x}^{k-1}\| + \|\bm{x}^{k-1}-\bm{x}^{k-2}\|\right)
\end{eqnarray}
for any $k \geq K$.
\end{lemma}
\beginproof
First, from \eqref{extrastep} and the first-order optimality condition for \eqref{subpro}, we have
\begin{numcases}{}
\bm{y}^{k-1} = \bm{x}^{k-1} + \bar{\beta}_{k-1} (\bm{x}^{k-1} - \bm{x}^{k-2}),  \label{extrastep2}    \\
0 \in \nabla f(\bm{y}^{k-1}) + \bar{\mu}_{k-1}(\bm{x}^{k} - \bm{y}^{k-1}) + \partial P(\bm{x}^{k}).  \label{suboptcond2}
\end{numcases}

Next, we consider the subdifferential of $H_{\delta}(\bm{u}, \bm{v}, \mu)$ at the point $(\bm{x}^{k}, \bm{x}^{k-1}, \bar{\mu}_{k-1})$ for $k\geq0$. Looking at the partial subdifferential with respect to $\bm{u}$, we have
\begin{eqnarray*}
\begin{aligned}
\partial_{\bm{u}} H_{\delta}(\bm{x}^{k}, \bm{x}^{k-1}, \bar{\mu}_{k-1})
&= \partial F(\bm{x}^{k}) + \frac{\delta\bar{\mu}_{k-1}}{2}(\bm{x}^{k} - \bm{x}^{k-1})  \\
&= \nabla f(\bm{x}^{k}) + \partial P(\bm{x}^{k}) + \frac{\delta\bar{\mu}_{k-1}}{2}(\bm{x}^{k} - \bm{x}^{k-1}) \\
&= \nabla f(\bm{y}^{k-1}) + \partial P(\bm{x}^{k}) + \bar{\mu}_{k-1}(\bm{x}^{k} - \bm{y}^{k-1}) + \nabla f(\bm{x}^{k}) - \nabla f(\bm{y}^{k-1}) \\
&\quad~ + \frac{\delta\bar{\mu}_{k-1}}{2}(\bm{x}^{k} - \bm{x}^{k-1}) - \bar{\mu}_{k-1}(\bm{x}^{k} - \bm{y}^{k-1})  \\
&\ni \nabla f(\bm{x}^{k}) - \nabla f(\bm{y}^{k-1}) + \frac{\delta\bar{\mu}_{k-1}}{2}(\bm{x}^{k} - \bm{x}^{k-1}) - \bar{\mu}_{k-1}(\bm{x}^{k} - \bm{y}^{k-1}),
\end{aligned}
\end{eqnarray*}
where the inclusion follows from \eqref{suboptcond2}. Similarly, we have
\begin{eqnarray*}
\begin{aligned}
\partial_{\bm{v}} H_{\delta}(\bm{x}^{k}, \bm{x}^{k-1}, \bar{\mu}_{k-1}) ~&=~ - \frac{\delta\bar{\mu}_{k-1}}{2}(\bm{x}^{k} - \bm{x}^{k-1}), \\
\partial_{\mu} H_{\delta}(\bm{x}^{k}, \bm{x}^{k-1}, \bar{\mu}_{k-1}) ~&=~ \frac{\delta}{4}\|\bm{x}^{k} - \bm{x}^{k-1}\|^2.
\end{aligned}
\end{eqnarray*}
Using the above relations, \eqref{extrastep2}, the global Lipschitz continuity of $\nabla f$ and the boundednesses of $\{\bar{\mu}_k\}$ and $\{\bar{\beta}_k\}$, there exists $c' > 0$ such that
\begin{eqnarray}\label{upbdpat}
\mathrm{dist}\left((0,0,0),\,\partial H_{\delta}(\bm{x}^{k}, \bm{x}^{k-1}, \bar{\mu}_{k-1})\right) \leq c'\left(\|\bm{x}^{k}-\bm{x}^{k-1}\| + \|\bm{x}^{k}-\bm{x}^{k-1}\|^2 + \|\bm{x}^{k-1}-\bm{x}^{k-2}\|\right).
\end{eqnarray}
Since $\|\bm{x}^{k+1} - \bm{x}^k\| \rightarrow 0$ (see Theorem \ref{subconvergence}(iv)), then there exists an integer $K>0$ such that $\|\bm{x}^{k+1} - \bm{x}^k\| \leq 1$ and hence $\|\bm{x}^{k+1} - \bm{x}^k\|^2 \leq \|\bm{x}^{k+1} - \bm{x}^k\|$ whenever $k \geq K$. This together with \eqref{upbdpat} completes the proof.
\endproof

In the next proposition, we discuss the KL exponent of our potential function $H_{\delta}$ defined in \eqref{defpofun}. This result can be viewed as a generalization of \cite[Theorem 3.6]{lp2017calculus} and the arguments involved are similar to those for \cite[Theorem 3.6]{lp2017calculus}, except that we have one more variable $\mu$ in $H_{\delta}$. For self-containedness, we provide the proof here.

\begin{proposition}[\textbf{KL exponent of $H_{\delta}$}]\label{expfun}
Suppose that Assumption \ref{assumfun} holds and the objective $F$ in \eqref{model} has the KL property at $\bar{\bm{x}}\in\mathrm{dom}\,\partial F$ with an exponent $\theta \in [\frac{1}{2}, 1)$. Let $\delta\geq0$ and $\mu_{\min} > 0$ be given constants. Then, for any $\bar{\mu} \geq \mu_{\min}$, the potential function $H_{\delta}$ defined in \eqref{defpofun} has the KL property at $(\bar{\bm{x}}, \bar{\bm{x}}, \bar{\mu})$ with an exponent $\theta$.
\end{proposition}
\beginproof
For $\delta=0$, the statement holds trivially since $H_{\delta}(\bm{u}, \bm{v}, \mu) \equiv F(\bm{u})$ if $\delta=0$. Thus, we only need to consider $\delta > 0$ in the following. Since $F$ has the KL property at $\bar{\bm{x}}\in\mathrm{dom}\,\partial F$ with an exponent $\theta \in [\frac{1}{2}, 1)$, there exist $\tilde{a},\varepsilon,\nu>0$ such that for any $\bm{x}$ satisfying $\bm{x} \in \mathrm{dom}\,\partial F$, $\|\bm{x}-\bar{\bm{x}}\|\leq\varepsilon$ and $F(\bar{\bm{x}}) < F(\bm{x}) < F(\bar{\bm{x}}) + \nu$, it holds that
\begin{eqnarray}\label{KLineq1}
\mathrm{dist}^{\frac{1}{\theta}}(0,\,\partial F(\bm{x})) \geq \tilde{a} (F(\bm{x}) - F(\bar{\bm{x}})).
\end{eqnarray}
Without loss of generality, we assume that $\varepsilon < \min\{\mu_{\min}, \,\frac{1}{2}\}$.

Next, consider any $(\bm{u},\bm{v},\mu)$ satisfying $\bm{u}\in\mathrm{dom}\,\partial F$, $\|\bm{u}-\bar{\bm{x}}\|\leq\varepsilon$, $\|\bm{v}-\bar{\bm{x}}\|\leq\varepsilon$, $|\mu-\bar{\mu}|\leq\varepsilon$ and $H_{\delta}(\bar{\bm{x}}, \bar{\bm{x}}, \bar{\mu}) < H_{\delta}(\bm{u}, \bm{v}, \mu) < H_{\delta}(\bar{\bm{x}}, \bar{\bm{x}}, \bar{\mu}) + \nu$. Note that, for any such $(\bm{u}, \bm{v}, \mu)$, we have
\begin{eqnarray*}
F(\bm{u}) \leq H_{\delta}(\bm{u}, \bm{v}, \mu) < H_{\delta}(\bar{\bm{x}}, \bar{\bm{x}}, \bar{\mu}) + \nu = F(\bar{\bm{x}}) + \nu.
\end{eqnarray*}
Thus, \eqref{KLineq1} holds for these $\bm{u}$ (if $F(\bm{u}) \leq F(\bar{\bm{x}})$, then \eqref{KLineq1} holds trivially). Moreover, for any such $(\bm{u},\bm{v},\mu)$, we have
\begin{eqnarray*}
\begin{aligned}
&~~ \mathrm{dist}^{\frac{1}{\theta}}((0,0,0),\,\partial H_{\delta}(\bm{u}, \bm{v}, \mu))   \\
&\geq c_0 \left(\left\|\frac{\delta\mu}{2}(\bm{u}-\bm{v})\right\|^{\frac{1}{\theta}} + \left(\frac{\delta}{4}\|\bm{u}-\bm{v}\|^2\right)^{\frac{1}{\theta}}
+\inf\limits_{\bm{\xi}\in\partial F(\bm{u})} \left\|\bm{\xi} + \frac{\delta\mu}{2}(\bm{u}-\bm{v})\right\|^{\frac{1}{\theta}}\right)  \\
&\geq c_0 \left(\left\|\frac{\delta\mu}{2}(\bm{u}-\bm{v})\right\|^{\frac{1}{\theta}} + \inf\limits_{\bm{\xi}\in\partial F(\bm{u})} \left\|\bm{\xi} + \frac{\delta\mu}{2}(\bm{u}-\bm{v})\right\|^{\frac{1}{\theta}}\right)  \\
&\geq c_0 \left(\left\|\frac{\delta\mu}{2}(\bm{u}-\bm{v})\right\|^{\frac{1}{\theta}} + \inf\limits_{\bm{\xi}\in\partial F(\bm{u})} \left(b_1\|\bm{\xi}\|^{\frac{1}{\theta}} - b_2\left\|\frac{\delta\mu}{2}(\bm{u}-\bm{v})\right\|^{\frac{1}{\theta}}\right)\right)  \\
&= c_0 \left(\frac{2^{2-\frac{1}{\theta}}(1-b_2)(\delta\mu)^{\frac{1}{\theta}-1}}{\tilde{a}} \cdot \frac{\tilde{a}\delta\mu}{4}\|\bm{u}-\bm{v}\|^{\frac{1}{\theta}} + b_1 \inf\limits_{\bm{\xi}\in\partial F(\bm{u})}\|\bm{\xi}\|^{\frac{1}{\theta}}\right)  \\
&\geq c_0\min\left\{\frac{2^{2-\frac{1}{\theta}}(1-b_2)(\delta\mu)^{\frac{1}{\theta}-1}}{\tilde{a}},\,b_1\right\}\cdot
\left(\frac{\tilde{a}\delta\mu}{4}\|\bm{u}-\bm{v}\|^{\frac{1}{\theta}} + \inf\limits_{\bm{\xi}\in\partial F(\bm{u})} \|\bm{\xi}\|^{\frac{1}{\theta}}\right)  \\
&\stackrel{(\mathrm{i})}{\geq} c_0\min\left\{\frac{2^{2-\frac{1}{\theta}}(1-b_2)(\delta\mu)^{\frac{1}{\theta}-1}}{\tilde{a}},\,b_1\right\}\cdot \left(\frac{\tilde{a}\delta\mu}{4}\|\bm{u}-\bm{v}\|^{\frac{1}{\theta}} + \tilde{a} (F(\bm{u}) - F(\bar{\bm{x}}))\right)  \\
&= c_0 \min\left\{2^{2-\frac{1}{\theta}}(1-b_2)(\delta\mu)^{\frac{1}{\theta}-1},\,\tilde{a}b_1\right\}\cdot\left(\frac{\delta\mu}{4}\|\bm{u}-\bm{v}\|^{\frac{1}{\theta}} + F(\bm{u}) - F(\bar{\bm{x}})\right)  \\
&\stackrel{(\mathrm{ii})}{\geq} c_0 \min\left\{2^{2-\frac{1}{\theta}}(1-b_2)(\delta(\mu_{\min}-\varepsilon))^{\frac{1}{\theta}-1},\,\tilde{a}b_1\right\}\cdot
\left(\frac{\delta\mu}{4}\|\bm{u}-\bm{v}\|^2 + F(\bm{u}) - F(\bar{\bm{x}})\right) \\
&= c_1 \Big(H_{\delta}(\bm{u}, \bm{v}, \mu) - H_{\delta}(\bar{\bm{x}}, \bar{\bm{x}}, \bar{\mu})\Big),
\end{aligned}
\end{eqnarray*}
where the existence of $c_0 > 0$ in the first inequality follows from Lemma \ref{normineq1}; the third inequality follows from Lemma \ref{normineq2} applied to $\|\bm{\xi} + \frac{\delta\mu}{2}(\bm{u}-\bm{v})\|^{\frac{1}{\theta}}$ with $b_1>0$ and $0<b_2<1$; the inequality (i) follows from \eqref{KLineq1}; the inequality (ii) follows from $\mu \geq \bar{\mu}-\varepsilon \geq \mu_{\min} - \varepsilon > 0$, $1<\frac{1}{\theta}\leq2$ and the observation that
\begin{eqnarray*}
\|\bm{u} - \bm{v}\| \leq \|\bm{u} - \bar{\bm{x}}\| + \|\bm{v} - \bar{\bm{x}}\| \leq 2 \varepsilon < 1.
\end{eqnarray*}
This completes the proof.
\endproof

Now, we are ready to discuss the local convergence rate of the PGels with $N=0$ under the additional assumption on the KL exponent of $F$. The proof is similar to that of Theorem \ref{conratenpg} but makes use of our potential function.

\begin{theorem}
Suppose that Assumption \ref{assumfun} holds, $\delta \in [0, \,1)$ is a nonnegative constant and the objective $F$ in problem \eqref{model} is a KL function with an exponent $\theta\in[\frac{1}{2}, 1)$. Let $\{\bm{x}^k\}$ and $\{\bar{\mu}_k\}$ be the sequences generated by Algorithm \ref{alg_PGels} with $N=0$, and let $\zeta$ be given in Theorem \ref{subconvergence}(iii). Then, the following statements hold.
\begin{itemize}
\item[{\rm (i)}] If $\theta = \frac{1}{2}$, there exist $\rho\in(0,1)$ and $c_1>0$ such that $|F(\bm{x}^k) - \zeta| \leq c_1 \rho^k$ for all large $k$;
\item[{\rm (ii)}] If $\theta \in \left(\frac{1}{2}, 1\right)$, there exists $c_2>0$ such that $|F(\bm{x}^k) - \zeta| \leq c_2\left(k-1\right)^{-\frac{1}{2\theta-1}}$ for all large $k$.
\end{itemize}
\end{theorem}
\beginproof
We first recall from \eqref{lscond} with $N=0$ that
\begin{eqnarray}\label{lscondN0}
H_{\delta}(\bm{x}^{k+1}, \bm{x}^k, \bar{\mu}_k) - H_{\delta}(\bm{x}^{k}, \bm{x}^{k-1}, \bar{\mu}_{k-1}) \leq -\frac{c}{2}\|\bm{x}^{k+1}-\bm{x}^{k}\|^2 \leq 0
\end{eqnarray}
for any $k\geq0$. Thus, the sequence $\{H_{\delta}(\bm{x}^{k}, \bm{x}^{k-1}, \bar{\mu}_{k-1})\}_{k=0}^{\infty}$ is obviously non-increasing. This together with Theorem \ref{subconvergence}(iii) implies that
\begin{eqnarray*}
\lim\limits_{k\rightarrow\infty} H_{\delta}(\bm{x}^{k}, \bm{x}^{k-1}, \bar{\mu}_{k-1})=\zeta, \quad \mathrm{and} \quad H_{\delta}(\bm{x}^{k}, \bm{x}^{k-1}, \bar{\mu}_{k-1})\geq\zeta ~\mathrm{for}~\mathrm{all}~k\geq0.
\end{eqnarray*}
Moreover, we have
\begin{eqnarray}\label{absfvalbd}
\begin{aligned}
|F(\bm{x}^k) - \zeta|
&= \left|H_{\delta}(\bm{x}^{k}, \bm{x}^{k-1}, \bar{\mu}_{k-1}) - \zeta - \frac{\delta\bar{\mu}_{k-1}}{4}\|\bm{x}^{k}-\bm{x}^{k-1}\|^2\right|  \\
&\leq \left|H_{\delta}(\bm{x}^{k}, \bm{x}^{k-1}, \bar{\mu}_{k-1}) - \zeta \right| + \frac{\delta\bar{\mu}_{k-1}}{4}\|\bm{x}^{k}-\bm{x}^{k-1}\|^2  \\
&\leq H_{\delta}(\bm{x}^{k}, \bm{x}^{k-1}, \bar{\mu}_{k-1}) - \zeta + \frac{\delta\bar{\mu}_{k-1}}{2c} \left(H_{\delta}(\bm{x}^{k-1}, \bm{x}^{k-2}, \bar{\mu}_{k-2})-H_{\delta}(\bm{x}^{k}, \bm{x}^{k-1}, \bar{\mu}_{k-1})\right) \\
&\leq H_{\delta}(\bm{x}^{k-1}, \bm{x}^{k-2}, \bar{\mu}_{k-2}) - \zeta + \frac{\delta\bar{\mu}_{k-1}}{2c} \left(H_{\delta}(\bm{x}^{k-1}, \bm{x}^{k-2}, \bar{\mu}_{k-2}) - \zeta\right) \\
&\leq {\textstyle\left(1+\frac{\delta\mu_{\max}}{2c}\right)}\left(H_{\delta}(\bm{x}^{k-1}, \bm{x}^{k-2}, \bar{\mu}_{k-2})-\zeta \right),
\end{aligned}
\end{eqnarray}
where the first equality follows from the definition of $H_{\delta}$ in \eqref{defpofun}, the first inequality follows from the triangle inequality, the second inequality  follows from \eqref{lscondN0} and the last three inequalities follow because $\{H_{\delta}(\bm{x}^{k}, \bm{x}^{k-1}, \bar{\mu}_{k-1})\}_{k=0}^{\infty}$ is non-increasing, $H_{\delta}(\bm{x}^{k}, \bm{x}^{k-1}, \bar{\mu}_{k-1})\geq\zeta$ and $\bar{\mu}_k \leq \mu_{\max}$ for all $k\geq0$. We next consider two cases.

\textbf{Case 1.} In this case, we suppose that $H_{\delta}(\bm{x}^{K_0}, \bm{x}^{K_0-1}, \bar{\mu}_{K_0-1}) = \zeta$ for some $K_0 \geq 0$. Since the sequence $\{H_{\delta}(\bm{x}^{k}, \bm{x}^{k-1}, \bar{\mu}_{k-1})\}_{k=0}^{\infty}$ is non-increasing, we must have $H_{\delta}(\bm{x}^{k}, \bm{x}^{k-1}, \bar{\mu}_{k-1}) = \zeta$ for all $k \geq K_0$. This together with \eqref{absfvalbd} proves statement (i) and (ii).

\textbf{Case 2.} From now on, we consider the case where $H_{\delta}(\bm{x}^{k}, \bm{x}^{k-1}, \bar{\mu}_{k-1}) > \zeta$ for all $k \geq 0$. From Theorem \ref{subconvergence}(i), we see that $\{\bm{x}^{k}\}$ is bounded and hence must have at least one cluster point. Let $\Gamma$ denote the set of cluster points of $\{\bm{x}^{k}\}$. Then, it follows from Proposition \ref{subfunseq}(ii) and $\ell(k)=k$ that $F \equiv \zeta$ on $\Gamma$.

Next, we consider the sequence $\{(\bm{x}^{k}, \bm{x}^{k-1}, \bar{\mu}_{k-1})\}_{k=0}^{\infty}$. In view of the boundedness of $\{\bar{\mu}_k\}$ (since $\mu_{\min}\leq\bar{\mu}_k\leq\mu_{\max}$) and Theorem \ref{subconvergence}(iv) which says that $\|\bm{x}^{k+1}-\bm{x}^{k}\|\rightarrow0$, it is not hard to show that the set of cluster points of $\{(\bm{x}^{k}, \bm{x}^{k-1}, \bar{\mu}_{k-1})\}_{k=0}^{\infty}$ is \textit{contained} in
\begin{eqnarray*}
\Upsilon:=\{(\bm{x},\bm{x},\mu): \bm{x}\in\Gamma, \,\mu_{\min}\leq\mu\leq\mu_{\max}\},
\end{eqnarray*}
which is a compact subset of $\mathrm{dom}\,\partial H_{\delta}$. Moreover, for any $(\bar{\bm{x}}, \bar{\bm{x}}, \bar{\mu})\in\Upsilon$ (hence $\bar{\bm{x}}\in\Gamma$ and $\mu_{\min}\leq\bar{\mu}\leq\mu_{\max}$), we have $H_{\delta}(\bar{\bm{x}}, \bar{\bm{x}},\bar{\mu})=F(\bar{\bm{x}})=\zeta$. Since $(\bar{\bm{x}}, \bar{\bm{x}}, \bar{\mu})\in\Upsilon$ is arbitrary, we conclude that $H_{\delta}\equiv\zeta$ on $\Upsilon$. On the other hand, since $F$ is a KL function with an exponent $\theta\in[\frac{1}{2}, 1)$, then it follows from Proposition \ref{expfun} that $H_{\delta}$ is a KL function with an exponent $\theta\in[\frac{1}{2}, 1)$ on $\Upsilon$. Using these facts and Proposition \ref{uniKL}, there exist $\varepsilon, \nu>0$ such that
\begin{eqnarray*}
\varphi'\left(H_{\delta}(\bm{u},\bm{v},\mu)-\zeta\right)\mathrm{dist}\left( (0,0,0),
\,\partial H_{\delta}(\bm{u},\bm{v},\mu)\right) \geq 1, \quad \mathrm{where}~\varphi(s)=as^{1-\theta}~\mathrm{for}~\mathrm{some}~a>0,
\end{eqnarray*}
for all $(\bm{u},\bm{v},\mu)$ satisfying $\mathrm{dist}((\bm{u},\bm{v},\mu), \,\Upsilon) < \varepsilon$ and
$\zeta < H_{\delta}(\bm{u},\bm{v},\mu) < \zeta + \nu$. Since $\Upsilon$ contains all the cluster points of $\{(\bm{x}^{k}, \bm{x}^{k-1}, \bar{\mu}_{k-1})\}_{k=0}^{\infty}$, then we have
\begin{eqnarray*}
\lim\limits_{k\rightarrow\infty}\mathrm{dist}((\bm{x}^{k}, \bm{x}^{k-1}, \bar{\mu}_{k-1}),\,\Upsilon) = 0.
\end{eqnarray*}
This together with $\lim\limits_{k\rightarrow\infty} H_{\delta}(\bm{x}^{k}, \bm{x}^{k-1}, \bar{\mu}_{k-1})=\zeta$ implies that there exists an integer $K_1\geq0$ such that $\mathrm{dist}((\bm{x}^{k}, \bm{x}^{k-1}, \bar{\mu}_{k-1}),\,\Upsilon) < \varepsilon$ and $\zeta < H_{\delta}(\bm{x}^{k}, \bm{x}^{k-1}, \bar{\mu}_{k-1}) < \zeta + \nu$ whenever $k \geq K_1$. Thus, for any $k\geq K_1$, we have
\begin{eqnarray}\label{klineq}
\varphi'\left(H_{\delta}(\bm{x}^{k}, \bm{x}^{k-1}, \bar{\mu}_{k-1})-\zeta\right)\mathrm{dist}\left( (0,0,0),
\,\partial H_{\delta}(\bm{x}^{k}, \bm{x}^{k-1}, \bar{\mu}_{k-1})\right) \geq 1.
\end{eqnarray}
For notational simplicity, let $\Delta_{H_{\delta}}^{k}:=H_{\delta}(\bm{x}^{k}, \bm{x}^{k-1}, \bar{\mu}_{k-1})
-\zeta$ and $\Delta_{\bm{x}^{k}}:=\bm{x}^{k+1}-\bm{x}^k$. Since the sequence $\{H_{\delta}(\bm{x}^{k}, \bm{x}^{k-1}, \bar{\mu}_{k-1})\}_{k=0}^{\infty}$ is non-increasing, then $\Delta_{H_{\delta}}^{k}$ is non-increasing, $\Delta_{H_{\delta}}^{k}>0$ for $k\geq0$ and $\Delta_{H_{\delta}}^{k}\rightarrow0$. We also let $\widebar{K}_1:=\max\{K_1, K\}$, where $K$ is given in Lemma \ref{pabound}. Then, for all $k \geq \widebar{K}_1$, we have
\begin{eqnarray}\label{keyineq}
\begin{aligned}
1&\leq \varphi'\left(\Delta_{H_{\delta}}^{k}\right)\mathrm{dist}\left( (0,0,0),
\,\partial H_{\delta}(\bm{x}^{k}, \bm{x}^{k-1}, \bar{\mu}_{k-1})\right) \\
&\leq \frac{a}{1-\theta}\cdot\left(\Delta_{H_{\delta}}^{k}\right)^{-\theta}\cdot \tilde{c}\left(\|\Delta_{\bm{x}^{k-1}}\| + \|\Delta_{\bm{x}^{k-2}}\|\right) \\
&\leq \frac{a\tilde{c}}{1-\theta}\cdot\left(\Delta_{H_{\delta}}^{k}\right)^{-\theta}\cdot
\sqrt{2\left(\|\Delta_{\bm{x}^{k-1}}\|^2 + \|\Delta_{\bm{x}^{k-2}}\|^2\right)} \\
&\leq \frac{a\tilde{c}}{1-\theta}\cdot\left(\Delta_{H_{\delta}}^{k}\right)^{-\theta}\cdot
\sqrt{\frac{4}{c}\big{(}H_{\delta}(\bm{x}^{k-2}, x^{k-3}, \bar{\mu}_{k-3})-H_{\delta}(\bm{x}^{k}, \bm{x}^{k-1}, \bar{\mu}_{k-1})\big{)}} \\
&= a_1\cdot\left(\Delta_{H_{\delta}}^{k}\right)^{-\theta}\cdot
\sqrt{ \Delta_{H_{\delta}}^{k-2}-\Delta_{H_{\delta}}^{k} }
\end{aligned}
\end{eqnarray}
where $a_1:=\frac{2a\tilde{c}}{(1-\theta)\sqrt{c}}>0$, the first inequality follows from \eqref{klineq}, the second inequality follows from \eqref{distbd}, the third inequality follows from $p+q \leq \sqrt{2(p^2+q^2)}$ for $p, q \geq 0$ and the last inequality follows from \eqref{lscondN0}. In the following, we consider two cases.

\textbf{(i)} $\theta=\frac{1}{2}$. Since $\Delta_{H_{\delta}}^{k}\rightarrow0$, there exists $K_2\geq0$ such that $\Delta_{H_{\delta}}^{k} \leq 1$ for all $k \geq K_2$. Then, for all $k \geq \widebar{K}_2:=\max\{K_2,\,\widebar{K}_1\}$, we see from \eqref{keyineq} that
\begin{eqnarray*}
\Delta_{H_{\delta}}^{k}=\left(\Delta_{H_{\delta}}^{k}\right)^{2\theta} \leq a_1^2\left(\Delta_{H_{\delta}}^{k-2}-\Delta_{H_{\delta}}^{k}\right).
\end{eqnarray*}
This implies that, for all $k \geq \widebar{K}_2+1$,
\begin{eqnarray*}
\Delta_{H_{\delta}}^{k-1} \leq \gamma \,\Delta_{H_{\delta}}^{k-3}
\leq \cdots \leq \gamma^{\frac{k-\widebar{K}_2-2}{2}} \,\Delta_{H_{\delta}}^{\widebar{K}_2},
\end{eqnarray*}
where $\gamma:=a_1^2/(1+a_1^2) < 1$. Then, from \eqref{absfvalbd}, we have
\begin{eqnarray*}
|F(\bm{x}^{k})-\zeta|
\leq {\textstyle\left(1+\frac{\delta\mu_{\max}}{2c}\right)} \Delta_{H_{\delta}}^{k-1}
\leq {\textstyle\left(1+\frac{\delta\mu_{\max}}{2c}\right)} \gamma^{\frac{k-\widebar{K}_2-2}{2}} \,\Delta_{H_{\delta}}^{\widebar{K}_2}
= c_1 \rho^k,
\end{eqnarray*}
where $c_1:=\left(1+\frac{\delta\mu_{\max}}{2c}\right)\gamma^{-\frac{\widebar{K}_2+2}{2}} \,\Delta_{H_{\delta}}^{\widebar{K}_2}$ and $\rho:=\sqrt{\gamma}<1$. This proves statement (ii).

\textbf{(ii)} $\frac{1}{2}<\theta<1$. We define $g(s):=s^{-2\theta}$ for $s\in(0,\,\infty)$. It is easy to see that $g$ is non-increasing. Then, for any $k \geq \widebar{K}_1$, we further consider the following two cases.
\begin{itemize}
\item If $g(\Delta_{H_{\delta}}^{k})\leq 2g(\Delta_{H_{\delta}}^{k-2})$, it follows from \eqref{keyineq} that
      \begin{eqnarray*}
      \begin{aligned}
      \frac{1}{a_1^2}
      &\leq (\Delta_{H_{\delta}}^{k})^{-2\theta}\cdot(\Delta_{H_{\delta}}^{k-2}-\Delta_{H_{\delta}}^{k})
       = g(\Delta_{H_{\delta}}^{k})\cdot (\Delta_{H_{\delta}}^{k-2}-\Delta_{H_{\delta}}^{k})
      \leq 2g(\Delta_{H_{\delta}}^{k-2})\cdot (\Delta_{H_{\delta}}^{k-2}-\Delta_{H_{\delta}}^{k})  \\
      &\leq 2\int^{\Delta_{H_{\delta}}^{k-2}}_{\Delta_{H_{\delta}}^{k}}\!\!\! g(s) \,\mathrm{d}s
       =\frac{2(\Delta_{H_{\delta}}^{k-2})^{1-2\theta}-2(\Delta_{H_{\delta}}^{k})^{1-2\theta}}{1-2\theta},
      \end{aligned}
      \end{eqnarray*}
      which, together with $1-2\theta<0$, implies that
      \begin{eqnarray}\label{difflbd1}
      (\Delta_{H_{\delta}}^{k})^{1-2\theta}-(\Delta_{H_{\delta}}^{k-2})^{1-2\theta} \geq \frac{2\theta-1}{2a_1^2}.
      \end{eqnarray}

\item If $g(\Delta_{H_{\delta}}^{k})\geq 2g(\Delta_{H_{\delta}}^{k-2})$, it is not hard to see that $(\Delta_{H_{\delta}}^{k})^{1-2\theta}\geq2^{\frac{2\theta-1}{2\theta}}(\Delta_{H_{\delta}}^{k-2})^{1-2\theta}$. Then, we have
      \begin{eqnarray}\label{difflbd2}
      (\Delta_{H_{\delta}}^{k})^{1-2\theta}-(\Delta_{H_{\delta}}^{k-2})^{1-2\theta}
      \geq \left(2^{\frac{2\theta-1}{2\theta}}-1\right)(\Delta_{H_{\delta}}^{k-2})^{1-2\theta}
      \geq \left(2^{\frac{2\theta-1}{2\theta}}-1\right)(\Delta_{H_{\delta}}^{\widebar{K}_1-2})^{1-2\theta},
      \end{eqnarray}
      where the last inequality follows from the facts that $\Delta_{H_{\delta}}^{k}$ is non-increasing and $1-2\theta<0$.
\end{itemize}

Thus, combining \eqref{difflbd1} and \eqref{difflbd2}, we obtain
\begin{eqnarray*}
(\Delta_{H_{\delta}}^{k})^{1-2\theta}-(\Delta_{H_{\delta}}^{k-2})^{1-2\theta}
\geq a_2 := \min\left\{\frac{2\theta-1}{2a_1^2},\,
\left(2^{\frac{2\theta-1}{2\theta}}-1\right)(\Delta_{H_{\delta}}^{\widebar{K}_1-2})^{1-2\theta}\right\}.
\end{eqnarray*}
Let $\pi_k = (k-\widebar{K}_1) \,\mathrm{mod}\,2$ for any $k\geq\widebar{K}_1$. Then, we have
\begin{eqnarray*}
\begin{aligned}
(\Delta_{H_{\delta}}^{k})^{1-2\theta}
&\geq (\Delta_{H_{\delta}}^{k})^{1-2\theta}-(\Delta_{H_{\delta}}^{\widebar{K}_1+\pi_k})^{1-2\theta}
=\sum^{k}_{j = \widebar{K}_1+\pi_k} \left((\Delta_{H_{\delta}}^{j})^{1-2\theta}-(\Delta_{H_{\delta}}^{j-2})^{1-2\theta}\right)\\
&\geq \frac{(k-\widebar{K}_1-\pi_k)a_2}{2} \geq \frac{a_2}{4}\,k,
\end{aligned}
\end{eqnarray*}
where the last inequality holds whenever $k \geq 2(\widebar{K}_1+1) \geq 2(\widebar{K}_1+\pi_k)$. Finally, using this relation and \eqref{absfvalbd}, we see that, for all $k \geq 2(\widebar{K}_1+1)+1$,
\begin{eqnarray*}
|F(\bm{x}^{k}) - \zeta|
\leq {\textstyle\left(1+\frac{\delta\mu_{\max}}{2c}\right)} \Delta_{H_{\delta}}^{k-1}
\leq {\textstyle\left(1+\frac{\delta\mu_{\max}}{2c}\right) \left(\frac{4}{a_2}\right)^{\frac{1}{2\theta-1}}}\,\left(k-1\right)^{-\frac{1}{2\theta-1}}
= c_2 \left(k-1\right)^{-\frac{1}{2\theta-1}},
\end{eqnarray*}
where $c_2:=\left(1+\frac{\delta\mu_{\max}}{2c}\right)\left(\frac{4}{a_2}\right)^{\frac{1}{2\theta-1}}$. This proves statement (ii).
\endproof

\section{Numerical experiments}\label{secnum}

In this section, we conduct some preliminary numerical experiments to test our PGels for solving the $\ell_1$ regularized logistic regression problem and the $\ell_{1\text{-}2}$ regularized least squares problem. All experiments are run in {\sc Matlab} R2016a on a 64-bit Laptop with an Intel Core i7-5600U CPU (2.60 GHz) and 8 GB of RAM equipped with Windows 10 OS.

\subsection{$\ell_1$ regularized logistic regression problem}\label{l1lognum}

In this subsection, we consider the $\ell_1$ regularized logistic regression problem
\begin{eqnarray}\label{l1logmodel}
\min\limits_{\tilde{\bm{x}}\in\mathbb{R}^n, x_0\in\mathbb{R}}~
F_{\log}(\bm{x}) := \sum^{m}_{i=1} \log(1+\exp(-b_i(\bm{a}_i^{\top}\tilde{\bm{x}}+x_0))) + \lambda\|\tilde{\bm{x}}\|_1,
\end{eqnarray}
where $\bm{x}:=(\tilde{\bm{x}}^{\top}, x_0)^{\top}\in\mathbb{R}^{n+1}$, $\bm{a}_i\in\mathbb{R}^n$, $b_i\in\{-1,1\}$ for $i=1,\cdots,m$ with $m < n$, and $\lambda>0$ is the regularization parameter. We further assume that $b_1, \cdots, b_m$ are not all the same. Let $C\in\mathbb{R}^{m\times (n+1)}$ be the matrix whose $i$-th row is given by $(\bm{a}_i^{\top},\,1)$. Then, we can rewrite \eqref{l1logmodel} in the form of \eqref{model} with
\begin{eqnarray*}
f(\bm{x}) = {\textstyle\sum^{m}_{i=1}}\log(1+\exp(-b_i(C\bm{x})_i)) \quad \mathrm{and} \quad P(\bm{x}) = \lambda\|\tilde{\bm{x}}\|_1.
\end{eqnarray*}
Moreover, one can check that $\nabla f$ is Lipschitz continuous with $L_f=0.25\|C\|^2$.

To apply our PGels, we also need to show that $F_{\log}$ is level-bounded\footnote{The level-boundedness of $F_{\log}$ was also shown using different way in the early version of \cite[Section 4.1]{wcp2017linear}, which is available at \url{https://arxiv.org/pdf/1512.09302v1.pdf}.}, i.e., $\mathrm{lev}_{\leq\alpha}F_{\log}:=\{\bm{x} : F_{\log}(\bm{x}) \leq \alpha\}$ is bounded (possibly empty) for every $\alpha\in\mathbb{R}$. Since $F_{\log}$ is nonnegative, then we only need to consider $\alpha\geq0$. For any $\bm{x}\in\mathrm{lev}_{\leq\alpha}F_{\log}$ with $\alpha\geq0$, due to the nonnegativities of $f$ and $P$, we have
\begin{numcases}{}
f(\bm{x})={\textstyle\sum^{m}_{i=1}} \log(1+\exp(-b_i(\bm{a}_i^{\top}\tilde{\bm{x}}+x_0))) \leq \alpha,  \label{fbdineq} \\
P(\bm{x})=\lambda\|\tilde{\bm{x}}\|_1 \leq \alpha.  \label{gbdineq}
\end{numcases}
Now, we define the index sets $\mathcal{I} := \{i : b_i = 1\}$ and $\mathcal{I}^c := \{i : b_i \neq 1\}$. Since $b_1, \cdots, b_m$ are not all the same, then both $\mathcal{I}$ and $\mathcal{I}^c$ are non-empty. Moreover, let $M:=\max\{\|\bm{a}_1\|_{\infty}$, $\cdots, \|\bm{a}_m\|_{\infty}\}$. Then, for any $i=1,\cdots,m$, it holds that
\begin{eqnarray*}
-\alpha M/\lambda \leq -\|\bm{a}_i\|_{\infty}\|\tilde{\bm{x}}\|_1 \leq \bm{a}_i^{\top}\tilde{\bm{x}}\leq\|\bm{a}_i\|_{\infty}\|\tilde{\bm{x}}\|_1\leq\alpha M/\lambda,
\end{eqnarray*}
where the first and last inequalities follow from \eqref{gbdineq}. Using the above relation, we further have
\begin{eqnarray}\label{ineqlog}
\begin{aligned}
\log(1+\exp(-\bm{a}_i^{\top}\tilde{\bm{x}}-x_0)) &\geq \log(1+\exp(-\alpha M/\lambda-x_0)), \quad i\in\mathcal{I}, \\
\log(1+\exp(\bm{a}_i^{\top}\tilde{\bm{x}}+x_0)) &\geq \log(1+\exp(-\alpha M/\lambda+x_0)), \quad i \in \mathcal{I}^c.
\end{aligned}
\end{eqnarray}
Thus, we see that
\begin{eqnarray*}
\begin{aligned}
&~~\log \big{(}(1+\exp(-\alpha M/\lambda-x_0))(1+\exp(-\alpha M/\lambda+x_0))\big{)}  \\
&=\log(1+\exp(-\alpha M/\lambda-x_0)) + \log(1+\exp(-\alpha M/\lambda+x_0)) \\
&\leq {\textstyle\sum_{i\in\mathcal{I}}} \log(1+\exp(-\alpha M/\lambda-x_0)) + {\textstyle\sum_{i\in\mathcal{I}^c}}
\log(1+\exp(-\alpha M/\lambda+x_0)) \\
&\leq {\textstyle\sum_{i\in\mathcal{I}}} \log(1+\exp(-\bm{a}_i^{\top}\tilde{\bm{x}}-x_0)) + {\textstyle\sum_{i\in\mathcal{I}^c}} \log(1+\exp(\bm{a}_i^{\top}\tilde{\bm{x}}+x_0)) \\
&=f(\bm{x}) \leq \alpha,
\end{aligned}
\end{eqnarray*}
where the first inequality follows because both $\mathcal{I}$ and $\mathcal{I}^c$ are non-empty, the second inequality follows from \eqref{ineqlog} and the last inequality follows from \eqref{fbdineq}. The above relation further implies that
\begin{eqnarray}\label{addineq1}
e^{\alpha}
\geq (1+\exp(-\alpha M/\lambda-x_0))(1+\exp(-\alpha M/\lambda+x_0))
= 1 + e^{-2\alpha M/\lambda} + e^{-\alpha M/\lambda} (e^{-x_0} + e^{x_0}).
\end{eqnarray}
Next, we consider the following two cases.
\begin{itemize}
\item $\alpha=0$. In this case, we see from \eqref{addineq1} that $e^{-x_0} + e^{x_0} \leq -1$, which cannot hold for any $x_0$. Thus, $\mathrm{lev}_{\leq0}F_{\log}$ is empty.

\item $\alpha>0$. In this case, it follows from \eqref{addineq1} that
      \begin{eqnarray*}
      e^{-x_0} + e^{x_0} \leq \widebar{M}:=e^{\alpha M/\lambda} (e^{\alpha} - e^{-2\alpha M/\lambda} - 1) \quad \Longrightarrow \quad |x_0| \leq \log \widebar{M}.
      \end{eqnarray*}
      This together with \eqref{gbdineq} shows that $\|\bm{x}\|$ is bounded and hence $\mathrm{lev}_{\leq\alpha}F_{\log}$ is bounded.
\end{itemize}
From the above, we see that $F_{\log}$ is level-bounded and hence our PGels is applicable.

\begin{figure}[ht]
\centering
\subfigure[$\mathrm{T}^{\max}=5$]{\includegraphics[width=5cm]{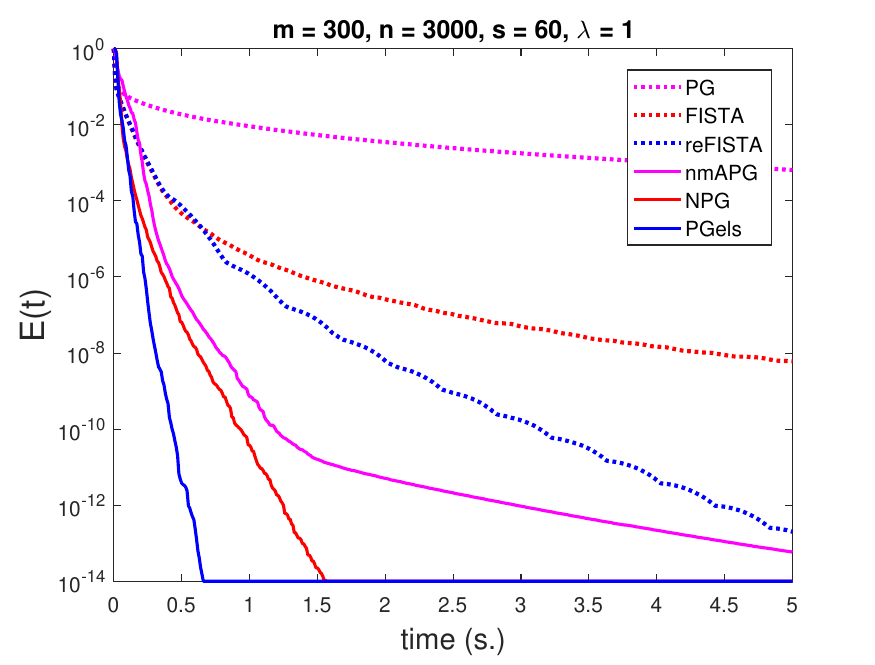}}
\subfigure[$\mathrm{T}^{\max}=15$]{\includegraphics[width=5cm]{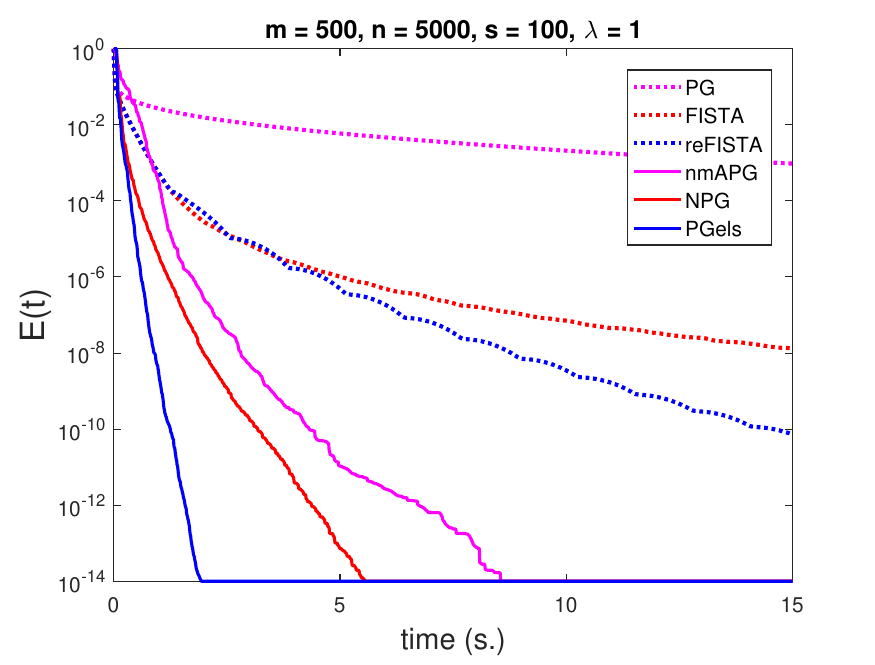}}
\subfigure[$\mathrm{T}^{\max}=50$]{\includegraphics[width=5cm]{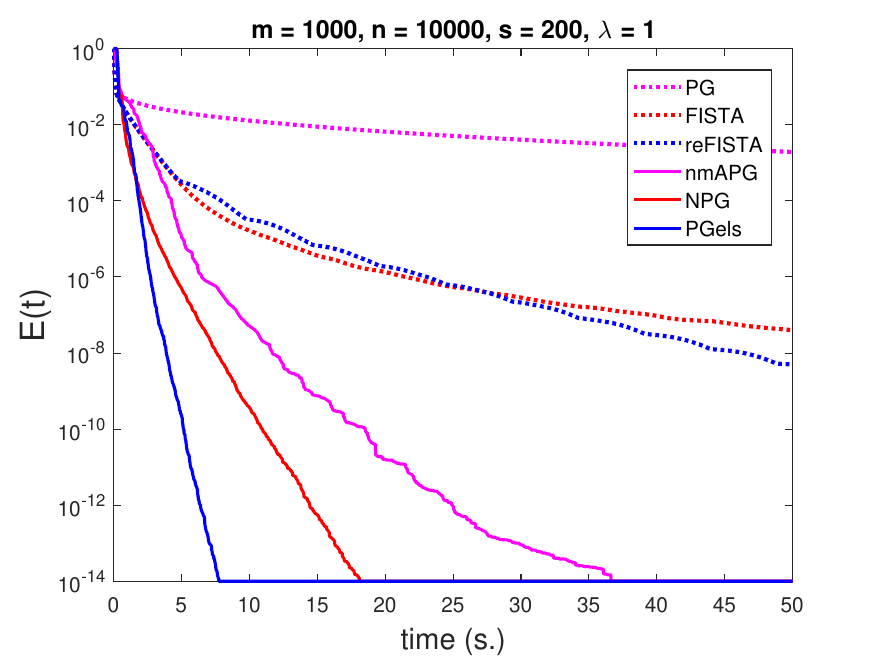}} \\
\subfigure[$\mathrm{T}^{\max}=10$]{\includegraphics[width=5cm]{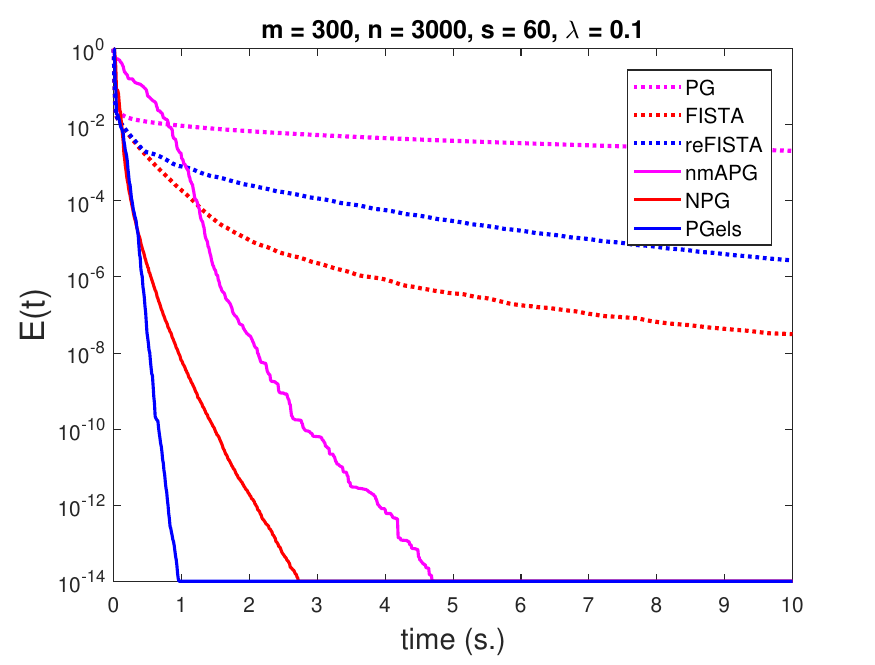}}
\subfigure[$\mathrm{T}^{\max}=30$]{\includegraphics[width=5cm]{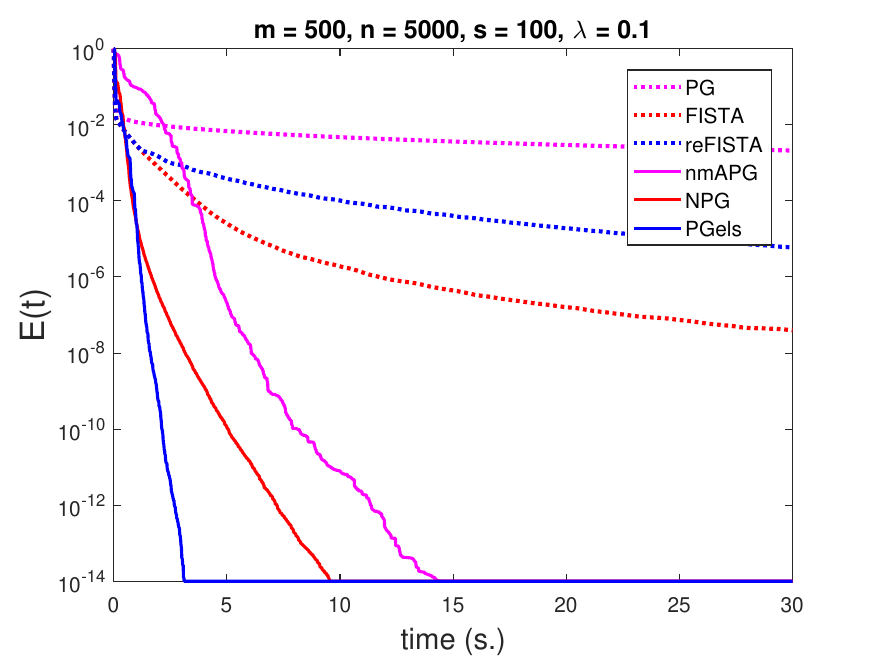}}
\subfigure[$\mathrm{T}^{\max}=100$]{\includegraphics[width=5cm]{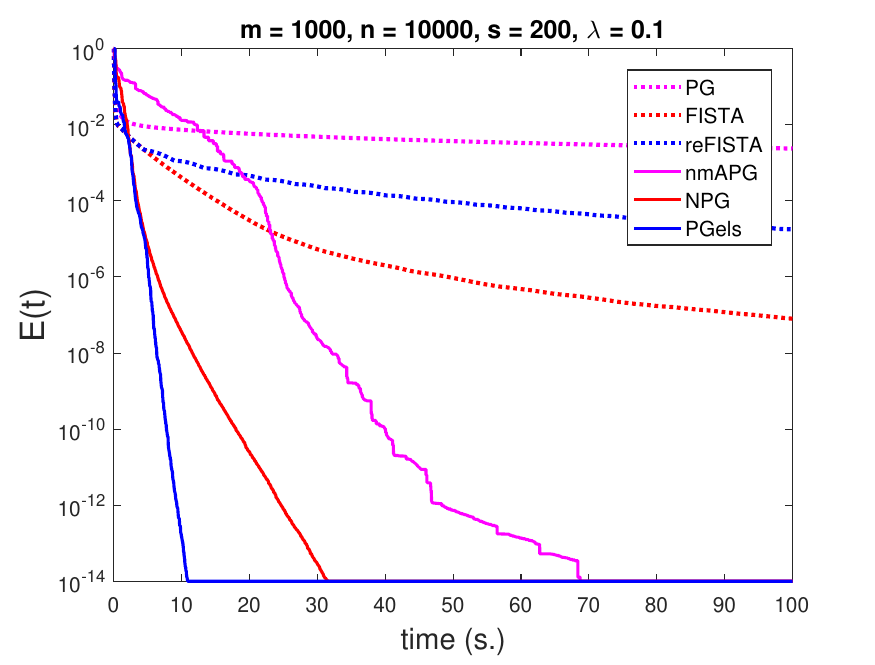}}
\caption{Average $E(t)$ of 10 independent trials of different algorithms for solving problem \eqref{l1logmodel}.}\label{Et_log}
\end{figure}

In our experiments, we will evaluate the PGels with $\delta=0.1$ (denoted by PGels) and the PGels with $\delta=0$ (denoted by NPG). For the PGels, we choose $\{\beta^0_k\}$ by \eqref{exparupdate} with $\beta^0_k$ in place of $\beta_k$. Moreover, for both PGels and NPG, we set $c=10^{-4}$, $\tau=2$, $\eta=0.8$, $N=2$, $\beta_{\max}=10$, $\mu_{\min}=10^{-6}$, $\mu_{\max}=\frac{L_f+2c}{1-\delta}$, $\mu_0^0=1$, and
\begin{eqnarray*}
\mu^0_k = \min\left\{\max\left\{\max\left\{ \frac{\langle \bm{y}^k-\bm{y}^{k-1}, \,\nabla f(\bm{y}^k) - \nabla f(\bm{y}^{k-1}) \rangle}{\|\bm{y}^k - \bm{y}^{k-1}\|^2}, \,0.5\bar{\mu}_{k-1}\right\}, \,\mu_{\min}\right\}, \,\mu_{\max}\right\}
\end{eqnarray*}
for $k\geq1$. We also compare PGels and NPG with PG, FISTA, FISTA with restart (reFISTA; see, e.g., \cite{bcg2011templates,oc2015adaptive,wcp2017linear}), and a non-monotone APG (nmAPG)\footnote{The implementations of nmAPG in our experiments are based on the original {\sc Matlab} codes, which are available at \url{http://www.cis.pku.edu.cn/faculty/vision/zlin/zlin.htm}} with line search \cite{ll2015accelerated}. For ease of future reference, we recall that FISTA for solving \eqref{l1logmodel} is given by
\begin{eqnarray}\label{FISTAscheme}
\left\{\begin{aligned}
&\beta_k = (t_{k-1} - 1)/t_k,  \\
&\bm{y}^k = \bm{x}^k + \beta_k (\bm{x}^k - \bm{x}^{k-1}),  \\
&\bm{x}^{k+1} = \mathrm{Prox}_{\frac{1}{L_f}P} \left( \bm{y}^k - \frac{1}{L_f}\nabla f(\bm{y}^k) \right),  \\
&t_{k+1} = \left(1+\sqrt{1+4t_k^2}\right)/2,
\end{aligned}\right.
\end{eqnarray}
with $\bm{x}^{-1}=\bm{x}^0$ and $t_{-1}=t_{0}=1$. Then, PG is given by \eqref{FISTAscheme} with $\beta_k \equiv 0$ and reFISTA is given by \eqref{FISTAscheme} with resetting $t_{k}=t_{k+1}=1$ whenever $k\,\mathrm{mod}\,\Delta K = 0$ or $\langle \bm{y}^{k} - \bm{x}^{k+1}, \, \bm{x}^{k+1} - \bm{x}^{k} \rangle > 0$. Moreover, nmAPG is developed based on \eqref{FISTAscheme} with a special monitor; see more details in \cite{ll2015accelerated}. In our experiments, we choose $\Delta K = 200$ for reFISTA (this restart interval has been observed in \cite{bcg2011templates} to have best performances). In addition, we initialize all algorithms at the origin and set the maximum running time\footnote{The maximum running time used in Sections \ref{l1lognum} and \ref{l1l2sec} does not include the time for computing the Lipschitz constant $L_f$ of $\nabla f$.} to $\mathrm{T}^{\max}$ for all algorithms. The specific values of $\mathrm{T}^{\max}$ are given in Figure\,\ref{Et_log}.

In the following experiments, we choose $\lambda \in \{1, \,0.1\}$ and consider $(m, n, s) = (100j, \,1000j, \,20j)$ for $j \in \{3, 5, 10\}$. For each triple $(m, n, s)$, we follow \cite[Section 4.1]{wcp2017linear} to randomly generate a trial as follows. First, we generate a matrix $A \in \mathbb{R}^{m\times n}$ with i.i.d. standard Gaussian entries. We then choose a subset $\mathcal{S}\subset\{1, \cdots, n\}$ of size $s$ uniformly at random and generate an $s$-sparse vector $\hat{\bm{x}}\in\mathbb{R}^{n}$, which has i.i.d. standard Gaussian entries on $\mathcal{S}$ and zeros on $\mathcal{S}^c$. Finally, we generate the vector $\bm{b} \in \mathbb{R}^{m}$ by setting $\bm{b} = \mathrm{sign}(A\hat{\bm{x}}+\hat{\epsilon}\mathbf{1})$, where $\hat{\epsilon}$ is chosen uniformly at random from $[0, 1]$ and $\mathbf{1}=(1,\cdots,1)^{\top}\in\mathbb{R}^{m}$.

To evaluate the performances of different algorithms, we follow \cite{gg2012accelerated,ypc2017a} to use an evolution of objective values. To introduce this evolution, we first define
\begin{eqnarray}\label{defek}
e(k):=\frac{F_{\log}(\bm{x}^k) - F_{\log}^{\min}}{F_{\log}(\bm{x}^0) - F_{\log}^{\min}},
\end{eqnarray}
where $F_{\log}(\bm{x}^k)$ denotes the objective value at $\bm{x}^k$ obtained by \textit{an} algorithm and $F_{\log}^{\min}$ denotes the minimum of the terminating objective values obtained among \textit{all} algorithms in a trial generated as above. For an algorithm, let $\mathcal{T}(k)$ denote the \textit{total} computational time (from the beginning) when it obtains $\bm{x}^k$. One can see that $\mathcal{T}(0)=0$ and $\mathcal{T}(k)$ is non-decreasing with respect to $k$. We now define the evolution of objective values obtained by a particular algorithm with respect to time $t$ as follows:
\begin{eqnarray*}
E(t):=\min\left\{e(k):k\in\{i: \mathcal{T}(i) \leq t\}\right\}.
\end{eqnarray*}
Note that $0 \leq E(t)\leq1$ (since $0 \leq e(k)\leq1$ for all $k$) and $E(t)$ is non-increasing with respect to $t$. It can be considered as a normalized measure of the reduction of the function value with respect to time. Then, one can take the average of $E(t)$ over several independent trials, and plot the average $E(t)$ within time $t$ for a given algorithm.

Figure\,\ref{Et_log} shows the average $E(t)$ of 10 independent trials of different algorithms for solving problem \eqref{l1logmodel}. From this figure, one can see that our PGels performs best in most cases in the sense that it takes less time to return a lower objective value.

\subsection{$\ell_{1\text{-}2}$ regularized least squares problem}\label{l1l2sec}

In this subsection, we consider the $\ell_{1\text{-}2}$ regularized least squares problem \cite{ylhx2015Minimization}:
\begin{eqnarray}\label{modell12}
\min\limits_{\bm{x}\in\mathbb{R}^n}~ F_{1\text{-}2}(\bm{x}):=\frac{1}{2}\|A\bm{x}-\bm{b}\|^2 + \lambda( \|\bm{x}\|_1 - \|\bm{x}\| ),
\end{eqnarray}
where $A\in\mathbb{R}^{m\times n}$, $\bm{b}\in\mathbb{R}^m$ and $\lambda>0$ is the regularization parameter. Obviously, this problem takes the form of $\eqref{model}$ with $f(\bm{x})=\frac{1}{2}\|A\bm{x}-\bm{b}\|^2$ and $P(\bm{x})=\lambda( \|\bm{x}\|_1 - \|\bm{x}\| )$. Moreover, $\nabla f$ is Lipschitz continuous with $L_f=\|A\|^2$ and $P(\bm{x})$ is a difference-of-convex regularizer. We further assume that $A$ does not have zero columns so that $F_{1\text{-}2}$ is level-bound; see \cite[Example 4.1(b)]{lp2017further} and \cite[Lemma 3.1]{ylhx2015Minimization}. Thus, our PGels is applicable.

\begin{figure}[ht]
\centering
\subfigure[$\mathrm{T}^{\max}=1$]{\includegraphics[width=5cm]{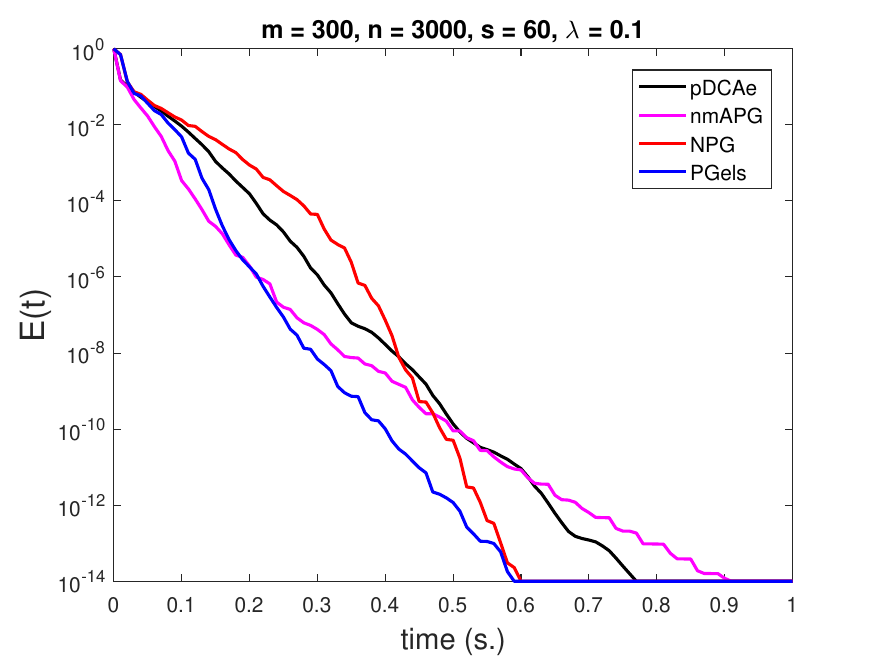}}
\subfigure[$\mathrm{T}^{\max}=3$]{\includegraphics[width=5cm]{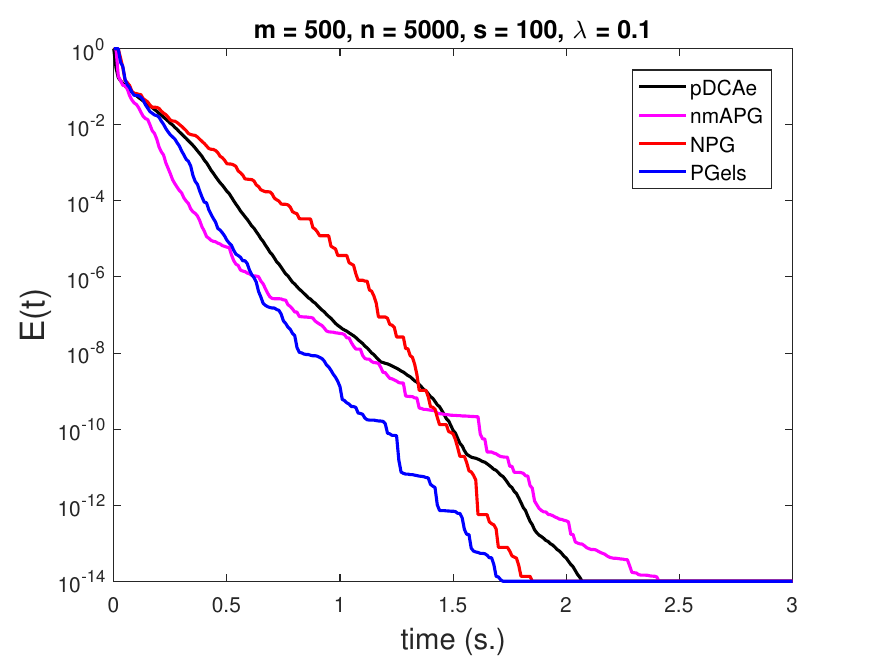}}
\subfigure[$\mathrm{T}^{\max}=15$]{\includegraphics[width=5cm]{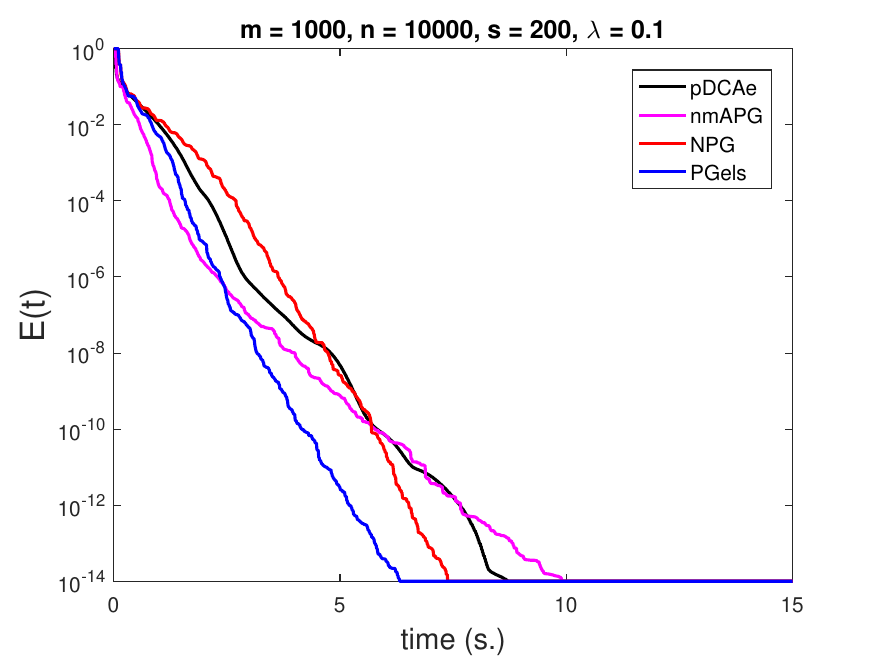}} \\
\subfigure[$\mathrm{T}^{\max}=15$]{\includegraphics[width=5cm]{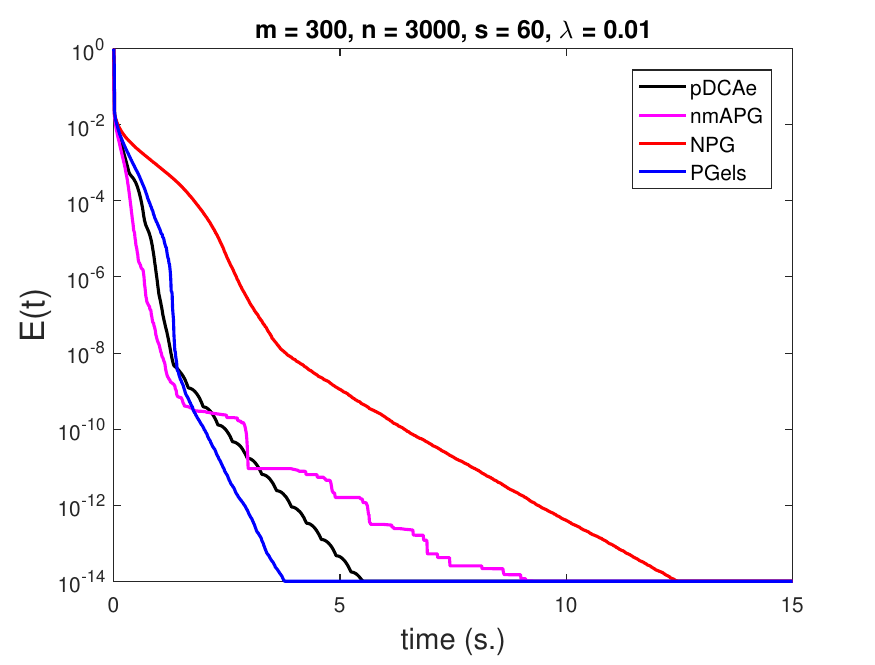}}
\subfigure[$\mathrm{T}^{\max}=50$]{\includegraphics[width=5cm]{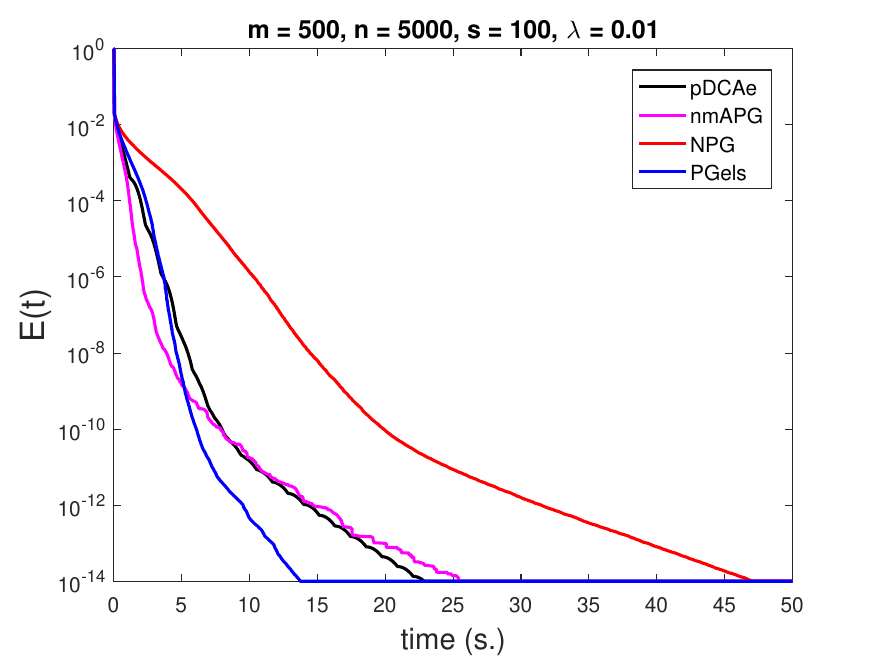}}
\subfigure[$\mathrm{T}^{\max}=150$]{\includegraphics[width=5cm]{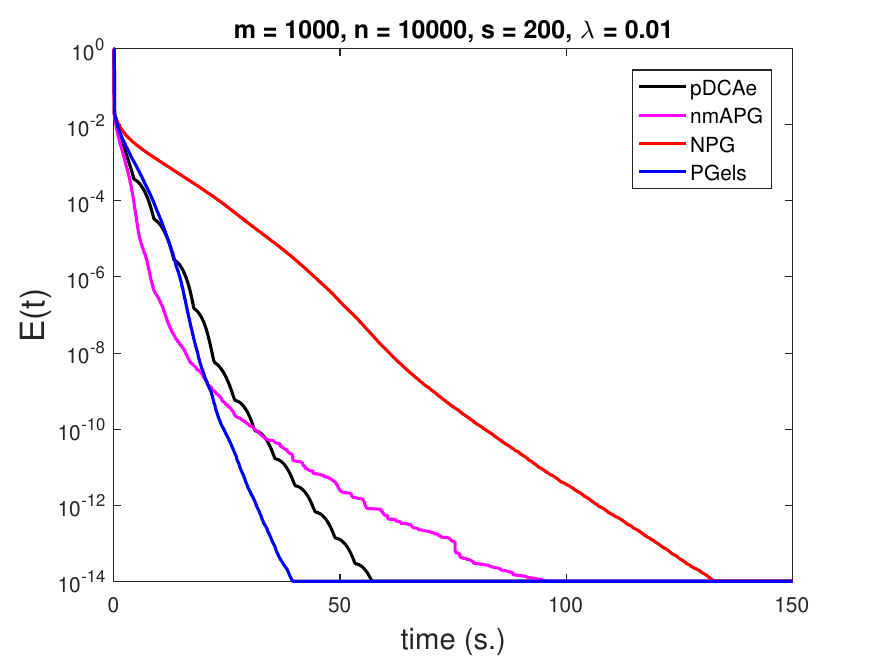}}
\caption{Average $E(t)$ of 10 independent trials of different algorithms for solving problem \eqref{modell12}.}\label{Et_l1l2}
\end{figure}

In this part of experiments, we compare four algorithms for solving \eqref{modell12}: PGels with $\delta=0.9$ (PGels), PGels with $\delta=0$ (NPG), nmAPG, and the proximal difference-of-convex algorithm with extrapolation ($\mathrm{pDCA}_{e}$)\footnote{The {\sc Matlab} codes of the $\mathrm{pDCA}_{e}$ for solving problem \eqref{modell12} are available at \url{http://www.mypolyuweb.hk/~tkpong/pDCAe_final_codes/}.} \cite{wcp2017a}. The $\mathrm{pDCA}_{e}$ for solving \eqref{modell12} is given as follows: choose proper extrapolation parameters $\{\beta_k\}$, let $\bm{x}^{-1} = \bm{x}^0$, and then at the $k$-th iteration,
\begin{eqnarray*}
\left\{\begin{aligned}
&\mathrm{Take}~\mathrm{any}~\bm{\xi}^k \in \lambda\partial \|\bm{x}^k\|~\mathrm{and}~\mathrm{compute}  \\
&\bm{y}^k = \bm{x}^k + \beta_k (\bm{x}^k - \bm{x}^{k-1}),  \\
&\bm{x}^{k+1} = \mathop{\mathrm{argmin}}\limits_{\bm{x}} \left\{\langle \nabla f(\bm{y}^k) - \bm{\xi}^k, \,\bm{x}\rangle + \frac{L_f}{2}\|\bm{x} - \bm{y}^k\|^2 + \lambda\|\bm{x}\|_1\right\}.
\end{aligned}\right.
\end{eqnarray*}
For PGels and NPG, we use the same parameter settings as in Section \ref{l1lognum}. For $\mathrm{pDCA}_{e}$, we follow \cite{wcp2017a} to choose $\{\beta_k\}$ as in reFISTA (see more details in Section \ref{l1lognum}). All algorithms are initialized at the origin and terminated by the maximum running time $\mathrm{T}^{\max}$. The specific values of $\mathrm{T}^{\max}$ are given in Fig.\,\ref{Et_l1l2}. In addition, as in Section \ref{l1lognum}, we also use the evolution of objective values (where $e(k)$ in \eqref{defek} is obtained by using $F_{1\text{-}2}$ in place of $F_{\log}$) to evaluate the performances of different algorithms.

In the following experiments, we choose $\lambda \in \{0.1, \,0.01\}$ and consider $(m, n, s) = (100j, \,1000j$, $20j)$ for $j \in \{3, \,5, \,10\}$. For each triple $(m, n, s)$, we follow \cite[Section 5]{wcp2017a} to randomly generate a trial as follows. We first generate a matrix $A \in \mathbb{R}^{m\times n}$ with i.i.d. standard Gaussian entries and then normalize $A$ so that the columns of $A$ have unit norms. We then uniformly at random choose a subset $\mathcal{S}$ of size $s$ from $\{1, \cdots, n\}$ and generate an $s$-sparse vector $\hat{\bm{x}}\in\mathbb{R}^{n}$, which has i.i.d. standard Gaussian entries on $\mathcal{S}$ and has zeros on $\mathcal{S}^c$. Finally, we set $\bm{b} = A\hat{\bm{x}}+0.01\cdot\hat{\bm{z}}$, where $\hat{\bm{z}}\in\mathbb{R}^{m}$ is a vector with i.i.d. standard Gaussian entries.

Figure\,\ref{Et_l1l2} shows the average $E(t)$ of 10 independent trials of different algorithms for solving \eqref{modell12}. From this figure, one can see that our PGels performs better than $\mathrm{pDCA}_{e}$ and is comparable with nmAPG. Note that $\mathrm{pDCA}_{e}$ is a difference-of-convex (DC) algorithm specifically designed for a class of DC problems taking the form of $\eqref{model}$, while our PGels can be applied for $\eqref{model}$ under more general scenarios. In addition, we observe that NPG performs worse in most cases. This situation was also observed in \cite{wcp2017a}. These observations, together with those observed in Section \ref{l1lognum}, show the potential advantage of combining extrapolation and non-monotone line search, which is the key motivation of this paper.

\section{Concluding remarks}\label{secconc}

In this paper, we considered a proximal gradient method with extrapolation and line search (PGels) for a composite optimization problem \eqref{model}, which is possibly nonconvex, nonsmooth, and non-Lipschitz. The basic idea of this method is to combine two simple and efficient acceleration techniques for PG, namely, extrapolation and non-monotone line search. We achieved this via the special potential function \eqref{defpofun}. By choosing proper parameters, PGels reduces to PG, PGe or NPG. We also established the global subsequential convergence for PGels. Specifically, under some mild conditions, we showed that the sequence generated by PGels is bounded and any cluster point of the sequence is a stationary point of \eqref{model}. In addition, by assuming that the objective in \eqref{model} is a Kurdyka-{\L}ojasiewicz function with an exponent $\theta$, we further studied the local convergence rate of two special cases of PGels, including NPG as one case. Finally, we conducted some numerical experiments to demonstrate the potential  advantage of combing two acceleration techniques.

%

%

\bibliographystyle{plain}
\bibliography{references/Ref_PGels}

\end{document}